\documentclass[11pt]{amsart} 
\usepackage{verbatim, latexsym, amssymb, amsmath,color}
\usepackage{epsfig}
\usepackage{tikz}
\usepackage{hyperref}
\usepackage{float}
\usetikzlibrary{matrix}
\usepackage{graphicx}
\def\R{\mathbb R}

\def\N{\mathbb N}

\newtheorem{thm}{Theorem}[section]

\newtheorem{prop}[thm]{Proposition}
\theoremstyle{remark}

\theoremstyle{definition}

\title{Min-max minimal hypersurfaces in non-compact manifolds}
\author{Rafael Montezuma}
\address{Instituto de Matem\'atica Pura e Aplicada (IMPA) \\ Estrada Dona Castorina 110 \\ 22460-320 Rio de Janeiro \\ Brazil}
\email{rafaelmc@impa.br}

\thanks{The author was partly supported by CNPq-Brazil and FAPERJ.}

\begin{document}

\begin{abstract}
{ {In this work we prove the existence of embedded closed minimal hypersurfaces in non-compact manifolds containing a bounded open subset with smooth and strictly mean-concave boundary and a natural behavior on the geometry at infinity. For doing this, we develop a modified min-max theory for the area functional following Almgren and Pitts' setting, to produce minimal surfaces with intersecting properties.}}
\end{abstract}
\maketitle
\setcounter{tocdepth}{1}

\section{Introduction}\label{introduction}

There is no immersed closed minimal surface in the Euclidean space $\R^3$. This fact illustrates the existence of simple geometric conditions creating obstructions for a Riemannian manifold to admit closed minimal surfaces. In the Euclidean space, we can see the obstruction coming in the following way: by the Jordan-Brouwer separation theorem every connected smooth closed surface $\Sigma^2 \subset \R^3$ divides $\R^3$ in two components, one of them bounded, which we denote $\Omega$. Start contracting a large Euclidean ball containing $\Omega$ until it touches $\Sigma$ the first time. Let $p \in \Sigma$ be a first contact point, then the maximum principle says that the mean curvature vector of $\Sigma$ at $p$ is non-zero and points inside $\Omega$. In particular, $\Sigma^2\subset \R^3$ is not minimal.    

In this work we consider two natural and purely geometric properties that imply that a complete non-compact Riemannian manifold $N$ admits a smooth closed embedded minimal hypersurface. In order to state the result, we introduce the following notation: we say that $N$ has the \textit{$\star_{k}$-condition} if there exists $p \in N$ and $R_0>0$, such that 
\begin{equation*}
\sup_{q \in B(p,R)} |\text{Sec}_N|(q) \leq R^k
\end{equation*}
and
\begin{equation*}
\inf_{q \in B(p,R)} inj_N(q) \geq R^{-\frac{k}{2}},
\end{equation*}
for every $R\geq R_0$, where $|\text{Sec}_N|(q)$ and $inj_N(q)$ denote, respectively, the maximum sectional curvature for $2$-planes contained in the tangent space $T_qN$ and the injectivity radius of $N$ at $q$. For instance, if $N$ has bounded geometry then the $\star_{k}$-condition holds for every positive $k$.

Our main result is:

\subsection*{Theorem A}\label{thm.A} \textit{ Let $(N^n,g)$ be a complete non-compact Riemannian manifold of dimension $n\leq 7$. Suppose:
\begin{itemize}
\item $N$ contains a bounded open subset $\Omega$, such that $\overline{\Omega}$ is a manifold with smooth and strictly mean-concave boundary;
\item $N$ satisfies the $\star_{k}$-condition, for some $k \leq \frac{2}{n-2}$.
\end{itemize}
Then, there exists a closed embedded minimal hypersurface $\Sigma^{n-1}\subset N$ that intersects $\Omega$.  
}

\begin{figure}[H]
  \centering
    \includegraphics[width=1.0\textwidth]{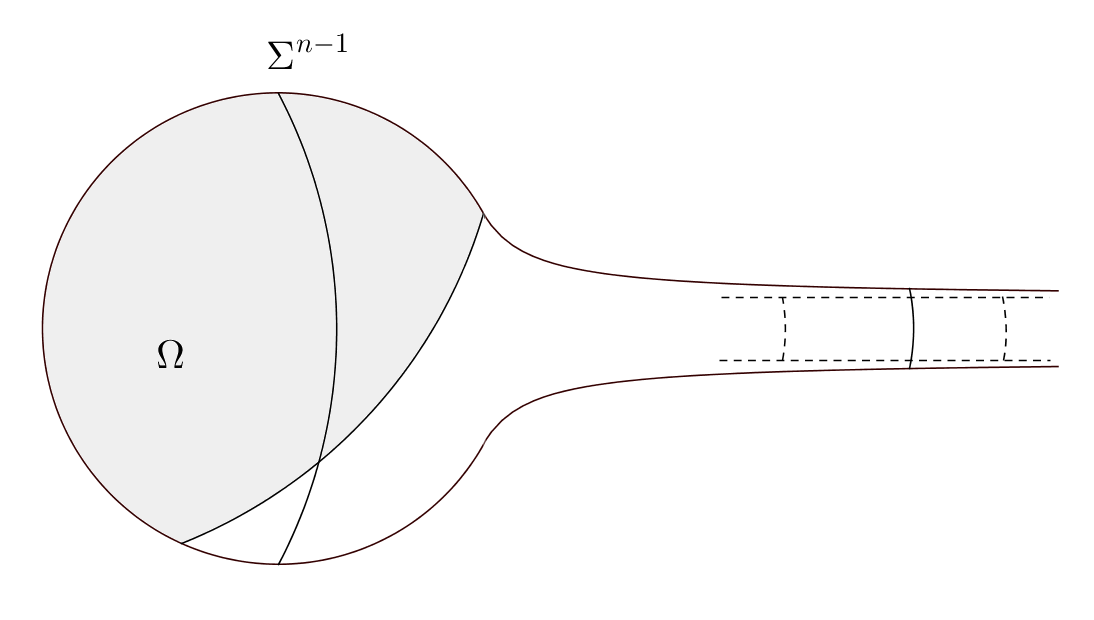}
\caption{A complete non-compact Riemannian manifold, asymptotic to a cylinder and containing a mean-concave open set $\Omega$. In this case, Theorem A could be applied.}\label{figure1}
\end{figure}

In the recent paper \cite{c-h-m-r}, Collin, Hauswirth, Mazet and Rosenberg prove that any complete non-compact hyperbolic three-dimensional manifold of finite volume admits a closed embedded minimal surface. These manifolds have a different behavior at infinity than those considered in Theorem A.

In our arguments, the geometric behavior of the ends of $M$ involving the $\star_{k}$-condition is used together with the monotonicity formula to provide a lower bound for the $(n-1)$-dimensional volume of minimal hypersurfaces in $M$ that, simultaneously, intersect $\Omega$ and contain points very far from it.

The hypothesis involving the mean-concave bounded domain $\Omega$ comes from the theory of closed geodesics in non-compact surfaces. In $1980$, Bangert proved the existence of infinitely many closed geodesics in a complete Riemannian surface $M$ of finite area and homeomorphic to either the plane, or the cylinder or the M\"{o}bius band, see \cite{bangert}. The first step in his argument is to prove that the finite area assumption implies the existence of locally convex neighborhoods of the ends of $M$. 

As a motivation for our approach, we briefly discuss how Bangert uses the Lusternik-Schnirelmann's theory in the case that $M$ is homeomorphic to the plane. If $C \subset M$ is a locally convex neighborhood of the infinity of $M$ whose boundary $\partial C \neq \varnothing$ is not totally geodesic, he proves that $M$ contains infinitely many closed geodesics intersecting $M - C$. To obtain one such curve, the idea is to apply the Lusternik-Schnirelmann's technique for a class $\Pi$ of paths $\beta$ defined on $[0,1]$ and taking values in a finite-dimensional subspace of the space of piecewise $C^1$ closed curves, with the properties that the curves $\beta_0$ and $\beta_1$ have image in the interior of $C$, being $\beta_0$ non-contractible in $C$ and $\beta_1$ contractible in $C$. The min-max invariant in this case is the number
\begin{equation*}
L(\Pi) = \inf_{\beta \in \Pi} \sup \{E(\beta_t) : \beta_{t}(S^1)\cap (M - \mathring{C}) \neq \varnothing\},
\end{equation*}
where $E(\gamma)$ denotes the energy of a map $\gamma : S^1\rightarrow M$, which is defined by
\begin{equation*}
E(\gamma) = \int_0^1 \vert \gamma^{\prime}(s)\vert^2 ds.
\end{equation*}
Then, achieve $L(\Pi)$ as the energy of a closed geodesic intersecting $M-C$.

To prove the Theorem A we develop a min-max method that is adequate to produce minimal hypersurfaces with intersecting properties. Let us briefly describe our technique. Let $(M^n,g)$ be a closed Riemannian manifold and $\Omega$ be an open subset of $M$. Consider a homotopy class $\Pi$ of one-parameter sweepouts of $M$ by codimension-one submanifolds. For each given sweepout $S = \{\Sigma_t\}_{t\in [0,1]} \in \Pi$, we consider the number
\begin{equation*}
L(S,\Omega) = \sup \{\mathcal{H}^{n-1}(\Sigma_t) : \Sigma_t \cap \overline{\Omega} \neq \varnothing\},
\end{equation*}
where $\mathcal{H}^{n-1}$ denotes the $(n-1)$-dimensional Hausdorff measure associated with the Riemannian metric. Define the width of $\Pi$ with respect to $\Omega$ to be
\begin{equation*}
L(\Pi,\Omega) = \inf \{L(S,\Omega) : S \in \Pi\}.
\end{equation*}

More precisely, our min-max technique is inspired by the discrete setting of Almgren and Pitts. The original method was introduced in \cite{alm2} and \cite{pitts} between the 1960's and 1980's, and has been used recently by Marques and Neves to answer deep questions in geometry, see \cite{marques-neves} and \cite{marques-neves-2}. The method consists of applications of variational techniques for the area functional. It is a powerful tool in the production of unstable minimal surfaces in closed manifolds. For instance, Marques and Neves, in the proof of the Willmore conjecture, proved that the Clifford Torus in the three-sphere is a min-max minimal surface. The min-max technique for the area functional appear also in a different setting, as introduced by Simon and Smith, in the unpublished work \cite{smith}, or Colding and De Lellis, in the survey paper \cite{c-dl}. Other recent developments on this theory can be found in \cite{dl-p}, \cite{dl-t}, \cite{ketover}, \cite{m.li}, \cite{xin}.

In the Almgren and Pitts' discrete setting, $\Pi$ is a homotopy class in $\pi^{\#}_1(\mathcal{Z}_{n-1}(M;\textbf{M}),\{0\})$. We define the width of $\Pi$ with respect to $\Omega$, $\textbf{L}(\Pi,\Omega)$, following the same principle as in the above discussion. Then, we prove:

\subsection*{Theorem B}
\textit{Let $(M^n,g)$ be a closed Riemannian manifold, $n\leq 7$, and $\Pi \in \pi^{\#}_1(\mathcal{Z}_{n-1}(M;\textbf{M}),\{0\})$ be a non-trivial homotopy class. Suppose that $M$ contains an open subset $\Omega$, such that $\overline{\Omega}$ is a manifold with smooth and strictly mean-concave boundary. There exists a stationary integral varifold $\Sigma$ whose support is a smooth embedded closed minimal hypersurface intersecting $\Omega$ and with $||\Sigma||(M) = \textbf{L}(\Pi,\Omega)$.
}

Consider the unit three-sphere $S^3 \subset \mathbb{R}^4$ and let $\Omega$ be a mean-concave subset of $S^3$. Assume that points in $\mathbb{R}^4$ have normal coordinates $(x,y,z,w)$. The min-max minimal surface $\Sigma$ produced from our method, starting with the homotopy class $\Pi$ of the standard sweepout $\{\Sigma_t\}$ of $S^3$, $\Sigma_t = \{w=t\}$ for $t \in [-1,1]$, is a great sphere. Indeed, it is obvious that 
$$\textbf{L}(\Pi,\Omega)\leq \max\{\mathcal{H}^2(\Sigma_t) : t \in [-1,1]\} = 4\pi.$$
This allows us to conclude that $\Sigma$ must be a great sphere, because these are the only minimal surfaces in $S^3$ with area less than or equal to $4\pi$.

\begin{figure}[H]
\centering
\includegraphics[width=1.0\textwidth]{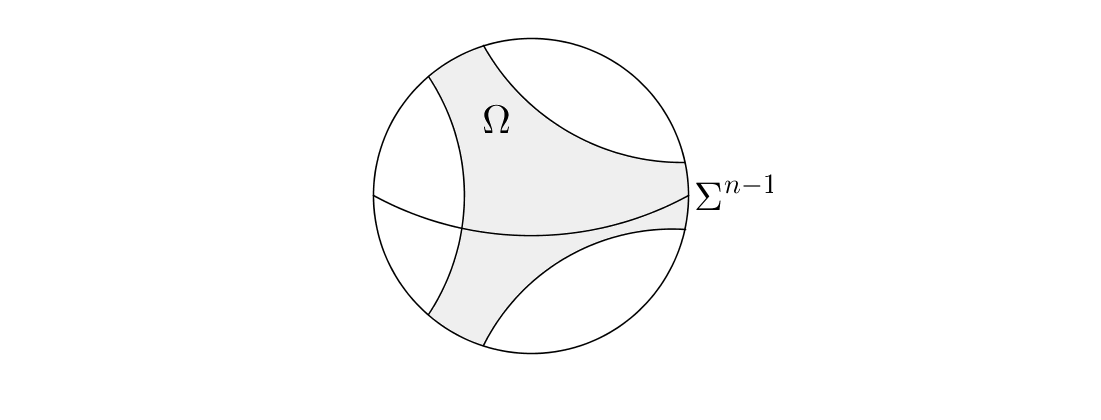}
\caption{Example of mean-concave domain $\Omega \subset S^3$ for which $\Sigma^{n-1}$ given by Theorem B is not entirely inside $\Omega$.}\label{figure2}
\end{figure}

The intersecting condition in Theorem B is optimal in the sense that it is possible that the support of $\Sigma$ is not entirely in $\overline{\Omega}$. We illustrate this with two examples of mean-concave subsets of the unit three-sphere $S^3 \subset \mathbb{R}^4$ containing no great sphere. The first example is the complement of three spherical geodesic balls, which can be seen in Figure \ref{figure2}. 
In order to introduce the second example, observe that $S^3 - (S^1\times \{(0,0)\}\cup \{(0,0)\}\times S^1)$ can be foliated by constant mean-curvature tori 
\begin{equation*}
S^1(t) \times S^1(\sqrt{1-t^2}) = \{(x,y,z,w) \in S^3 : x^2 + y^2 = t^2 \text{ and } z^2 + w^2 = 1-t^2\},
\end{equation*}
for $0<t<1$. Each such surface divides $S^3$ in two components
\begin{equation*}
\Omega_{+}(t) = S^3\cap \{x^2 + y^2 > t^2\} \quad \text{and} \quad \Omega_{-}(t) = S^3\cap \{ x^2 + y^2 < t^2\}.
\end{equation*}
We claim that neither $\Omega_{+}(t)$ nor $\Omega_{-}(t)$ contain great spheres. In the case of $\Omega_{+}(t)$, we only observe that the great sphere $S^3\cap \{x=0\}$ contains no great circle inside $\Omega_{+}(t)$. But any two great spheres in $S^3$ have a great circle in common, at least. The case of $\Omega_{-}(t)$ is analogous. A simple calculation yields that $S^1(t) \times S^1(\sqrt{1-t^2}) \subset S^3$ has mean-curvature vector given by
\begin{equation*}
\overrightarrow{H}(x,y,z,w) = (1-2t^2)\left(-\frac{x}{t^2}, -\frac{y}{t^2}, \frac{z}{1-t^2}, \frac{w}{1-t^2}\right).
\end{equation*} 
In particular, if $t= 1/\sqrt{2}$ the torus $S^1(1/\sqrt{2})\times S^1(1/\sqrt{2})$ is a minimal surface, known as the Clifford torus. Moreover, $\overrightarrow{H}$ points inside $\Omega_{-}(t)$, if $t< 1/\sqrt{2}$, and points inside $\Omega_{+}(t)$, if $t>1/\sqrt{2}$. Therefore, either $\Omega_{+}(t)$ or $\Omega_{-}(t)$ is strictly mean-concave and contains no great sphere.

In future work, we will adapt our methods to the set up of Simon and Smith. This different approach is important for a better understanding of the geometrical and topological behavior of the produced minimal surfaces.


\subsection*{Acknowledgments} The results contained in this paper are based partially on the author's Ph.D. thesis under the guidance of Professor Fernando Cod\'a Marques. It is a pleasure to show my gratefulness to him for his support during the preparation of this work.


\subsection*{Organization} The content of this paper is organized as follows:

In Section \ref{outline}, we outline the main ideas to prove Theorem B, avoiding the technical details and the language from the geometric measure theory.

In Section \ref{section.pre}, we recall some definitions from geometric measure theory, state the maximum principle for varifolds and introduce the discrete maps in the space of currents.

In Section \ref{min-max.section}, we develop the min-max theory for intersecting slices. We give the definitions, state the main results of our method and compare them with the results obtained by the classical setting of Almgren and Pitts. The proofs of the results stated in Section \ref{min-max.section} are done in Sections \ref{sec.interp.discrt.cont} to \ref{sect-proof.of.am}.

In Section \ref{sec.interp.discrt.cont}, we describe the interpolation results that we use. They are a powerful tool to deal with discrete maps in the space of currents.

In Section \ref{section_natural sweepouts}, we develop a deformation argument for boundaries of integral $n$-currents on compact $n$-dimensional manifolds with boundary. This is useful to prove the existence of critical sequence in our min-max setting.

In Sections \ref{small.mass.section} and \ref{sect.small.mass}, we prove a deformation lemma for discrete maps in the space of codimension-one integral currents. This result plays an important role in the subsequent sections.

In Section \ref{creating.discrete.sweepouts}, we describe how to create discrete sweepouts in the setting of Almgren and Pitts out of a continuous one.

In Section \ref{proof.pulltight}, we adapt the pull-tight procedure to our setting.

In Section \ref{sect-proof.of.am}, we prove the existence of intersecting and almost minimizing critical varifolds in a given homotopy class.

In Section \ref{exist.non-compact}, we apply our min-max method to prove the existence of closed minimal surfaces in complete non-compact manifolds, Theorem A.


\section{Main ideas of the min-max theory for intersecting slices}\label{outline}

In this section, we outline the proof of our min-max result, Theorem B. We present the main ideas, omitting the technical issues and the language and some ingredients from the geometric measure theory. 

Let $(M^n,g)$ be an orientable closed Riemannian manifold and $\Omega \subset M$ be a connected open subset with smooth and strictly mean-concave boundary. We consider sweepouts $S = \{\Sigma_t\}_{t\in [0,1]}$ of codimension-one submanifolds $\Sigma_t \subset M$, such that at times $t =0$ and $1$ the slices  $\Sigma_t$ degenerate to points in $M$. For example, one can consider the sweepout of the level sets $\Sigma_t = f^{-1}(t)$ of a Morse function $f : M \rightarrow [0,1]$. The key difference between our method  and the original one by Almgren and Pitts is that we see only the slices intersecting $\overline{\Omega}$. Use the following definition:
\begin{equation*}
\text{dmn}_{\Omega}(S) = \{t \in [0,1] : \Sigma_{t}\cap \overline{{\Omega}} \neq \varnothing\}.
\end{equation*}
There is a notion of homotopy between sweepouts of $M$. Once a homotopy class $\Pi$ of sweepouts is fixed, one can run the min-max. For each sweepout $S  = \{\Sigma_t\}_{t\in [0,1]} \in \Pi$, we consider
\begin{equation*}
\textbf{L}(S,\Omega) = \sup\{\mathcal{H}^{n-1}(\Sigma_t) : t \in \text{dmn}_{\Omega}(S)\}.
\end{equation*}
And we define the width of $\Pi$ with respect to $\Omega$ to be
\begin{equation*}
\textbf{L}(\Pi,\Omega) = \inf \{ \textbf{L}(S,\Omega): S \in \Pi\}.
\end{equation*}
The standard setting of Almgren and Pitts coincides with ours when $\Omega$ is the whole $M$. In this case, the min-max philosophy was to obtain sequences $S_k = \{\Sigma^k_t\}_{t\in [0,1]} \in \Pi$, for $k = 1,2, \ldots$, and $\{t_k\}_{k\in \mathbb{N}} \subset [0,1]$ such that
\begin{equation*}
\lim_{k\rightarrow \infty} \mathcal{H}^{n-1}(\Sigma^k_{t_k})  = \lim_{k\rightarrow \infty} \textbf{L}(S_k,M)= \textbf{L}(\Pi,M) = \inf_{\{\Sigma_t\} \in \Pi} \ \ \sup\quad \mathcal{H}^{n-1}(\Sigma_t)
\end{equation*}
and
\begin{equation*}
\Sigma^k_{t_k} \rightharpoonup V, \quad \text{in the sense of varifolds},
\end{equation*}
for some stationary integral varifold $V$, whose support is a union of closed embedded smooth minimal hypersurfaces. In our approach, we show that it is possible to attain the width $\textbf{L}(\Pi,\Omega)$ via a sequence of slices $\Sigma^k_{t_k}$ as above and with the extra properties that $t_k \in \text{dmn}_{\Omega}(S_k)$ and $\{\Sigma^k_{t_k}\}_{k\in \mathbb{N}}$ converges to a varifold $V$ whose support intersects $\Omega$.

Despite we are avoiding the language of the geometric measure theory in this overview section, it is important to observe that the min-max minimal hypersurface arises as a limit of big slices in the weak sense of varifolds. Also, these slices are not hypersurfaces in general, they are codimension-one integral currents in $M^n$, which can be regarded as Lipschitz oriented $(n-1)$-submanifolds, with integer multiplicity, with no boundary. On the other hand, a $(n-1)$-varifold $V$ in $M$ is a Radon measure on the Grassmannian of unoriented $(n-1)$-subspaces on $M$. For example, given an embedded hypersurface $\Sigma^{n-1}\subset M^n$  and a positive measurable function $\theta : \Sigma \rightarrow \mathbb{R}$, we can consider the linear map
\begin{equation*}
V(\varphi) := \int_{\Sigma} \varphi(x,T_x\Sigma)\theta(x) d\Sigma(x), \quad \text{where } \varphi \in C_0(G_{n-1}(M)),
\end{equation*}
where $C_0(G_{n-1}(M))$ denotes the space of the continuous functions with compact support over the Grassmannian of unoriented $(n-1)$-subspaces on $M$. By the Riesz representation theorem, $V$ can be thought as a $(n-1)$-varifold in $M$, the varifold associated to $\Sigma^{n-1}$ of density $\theta$. In particular, one can associate a varifold for each current simply by forgetting the orientation. Varifolds also have a natural notion of mass, which coincides with $\mathcal{H}^{n-1}$, the Hausdorff measure with respect to $g$, for hypersurfaces. The critical points for the mass functional are called stationary varifolds.  

The first thing that we check is the existence of homotopy classes with $\textbf{L}(\Pi,\Omega)>0$. The existence of non-trivial homotopy classes comes back to Almgren's work, see \cite{alm1}. He proved that the set of homotopy classes of sweepout is isomorphic to the homology group $H_n(M^n,\mathbb{Z})$. Moreover, it is possible to prove the existence of positive constants $\alpha_0 = \alpha_0(M)$ and $r_0 = r_0(M)$ with the property that: given $p\in M$, $0<r\leq r_0$ and a sweepout $S = \{\Sigma_t\}_{t\in [0,1]}$ in a non-trivial homotopy class, we have
\begin{equation*}
\sup_{t\in [0,1]} \mathcal{H}^{n-1}(\Sigma_t \cap B(p,r)) \geq \alpha_0 r^{n-1},
\end{equation*}
where $B(p,r)$ denotes the geodesic ball of $M$, with radius $r$ and centered at $p$. This implies that, if we choose $\overline{B(p,r)} \subset \Omega$, then
\begin{equation*}
\textbf{L}(S,\Omega) = \sup_{t \in \text{dmn}_{\Omega}(S)} \mathcal{H}^{n-1}(\Sigma_t) \geq \sup_{t\in [0,1]} \mathcal{H}^{n-1}(\Sigma_t \cap B(p,r)) \geq \alpha_0 r^{n-1}.
\end{equation*}
Therefore, $\textbf{L}(\Pi,\Omega)>0$, for non-trivial homotopy classes $\Pi$. The existence of such numbers $\alpha_0$ and $r_0$ is due to Gromov, see Section $4.2.B$ in \cite{gromov}. These lower bounds for the width also appear in \cite{guth} and Section $8$ of \cite{marques-neves-2}.

From now on we fix a non-trivial homotopy class $\Pi$. We say that a sequence $\{S_k\}_{k\in \mathbb{N}}\subset \Pi$ is critical with respect to $\Omega$ if we have
\begin{equation}
\textbf{L}(\Pi,\Omega) = \lim_{k\rightarrow \infty} \textbf{L}(S_k,\Omega)
\end{equation}
and
\begin{equation}
\sup\{\mathcal{H}^{n-1}(\Sigma^k_t) : k \in \mathbb{N} \text{ and } t \in [0,1]\}< \infty.
\end{equation}
We adopt the notation that $S_k = \{\Sigma^k_t\}_{t\in [0,1]}$. The second condition above is trivial in the classical case $\Omega = M$, because the measure of all slices are accounted on the width. To assure that it is possible to achieve the width keeping controlled the $\mathcal{H}^{n-1}$-measure of all slices, we develop a deformation process for boundaries of open sets in compact manifolds with boundary.

Given a critical sequence with respect to $\Omega$, $S = \{S_k\}_{k \in \mathbb{N}}$, we define its critical set $\mathcal{C}(S,\Omega)$ to be the set of limit varifolds $V = \lim_{k} \Sigma^k_{t_k}$, for some sequence of $t_k \in \text{dmn}_{\Omega}(S_k)$. Formally, this limit is in the weak sense of varifolds. Intuitively, we can use the compactness theorem for varifolds with uniformly bounded masses to pass to a limit sweepout, also denoted by $S = \{\Sigma_t\}_{t\in [0,1]} \in \Pi$. Observe that we have $\textbf{L}(S,\Omega) = \textbf{L}(\Pi,\Omega)$. In this case, we say that $S$ is an optimal sweepout in $\Pi$ with respect to $\Omega$ and, moreover, we can identify the critical set $\mathcal{C}(S,\Omega)$ of the critical sequence $S$ with the set of the maximal intersecting slices of the sweepout $S$. 

This optimal sweepout does not belong necessarily to the same category of sweepouts as the $S_k$'s, but for our purposes, let us suppose for a while that there exists such a $S^* = \{\overline{\Sigma}_t\}_{t\in [0,1]} \in \Pi$, optimal with respect to $\Omega$.

Next, we deform $S^*$ to obtain a better optimal sweepout $S = \{\Sigma_{t}\}_{t\in [0,1]} \in \Pi$ with the extra property that: if $t \in \text{dmn}_{\Omega}(S)$ and $\mathcal{H}^{n-1}(\Sigma_t) = \textbf{L}(S,\Omega)$, then either $\Sigma_t$ is stationary in $M$ or $\Sigma_t \cap \overline{\Omega} \subset \partial \Omega$. This is inspired by the original pull-tight argument, see Theorem $4.3$ in \cite{pitts} and Proposition $8.5$ in \cite{marques-neves}. More precisely, we construct continuous slice deformations $\{H(s,t)\}_{s,t \in [0,1]}$, starting with $H(0,t) = \overline{\Sigma}_t$. The key property being that $$\mathcal{H}^{n-1}(H(1,t)) < \mathcal{H}^{n-1}(t),$$ unless $\overline{\Sigma}_t$ is either stationary or $\text{spt}(\overline{\Sigma}_t)\cap \Omega = \varnothing$, for which $H(s,t)=H(0,t)$, for all $0\leq s\leq 1$. Then we consider the sweepout $S = \{\Sigma_t\}_{t \in [0,1]}$ given by $\Sigma_t = H(1,\overline{\Sigma}_t)$. The main improvement or our pull-tight deformation $H$ is that it keeps unmoved the non-intersecting slices, allowing us to conclude that $S$ is also optimal with respect to $\Omega$.

Recall that, in general, $S^*$ is not a true sweepout in the homotopy class $\Pi$ and we have to deal with a critical sequence $S^* = \{S^*_k\}_{k\in \mathbb{N}}\subset \Pi$. What we are able to obtain in terms of sequences is a new critical sequence with respect to $\Omega$, $S = \{S_k\}_{k\in \mathbb{N}} \subset \Pi$, for which each $V \in \mathcal{C}(S,\Omega)$ is either stationary in $M$ or its support does not intersect $\Omega$.

The final part is to prove that some $V \in \mathcal{C}(S,\Omega)$, for the above sequence $S\subset \Pi$, is at the same time intersecting and smooth. Since the intersecting critical $V \in \mathcal{C}(S,\Omega)$ are stationary, some component of $\text{spt}(V)$ is a closed embedded minimal hypersurface in $M$ and it intersects $\Omega$. In order to prove the smoothness assumption, we prove the following claim:

\subsection{Claim}\label{claim_outline}\textit{There exists an intersecting critical varifold $V$ that is almost minimizing in small annuli, in the sense introduced by Pitts.}

For the precise notation, see definition \ref{am.varifolds} and Theorem \ref{existence.almost.mini.varifolds} in Section \ref{min-max.section}, or Theorem $4.10$ in Pitts book \cite{pitts}. The almost minimizing varifolds are natural objects in the min-max theory, which can be arbitrarily approximated by locally mass-minimizing currents. In that same work \cite{pitts}, Pitts proved that stationary varifolds that are almost minimizing in annuli are smooth. This will conclude our argument.

The proof of Claim \ref{claim_outline} is by contradiction. Assuming that no $V \in \mathcal{C}(S,\Omega)$ satisfies our assumption, we deform the critical sequence $S \subset \Pi$ to obtain strictly better competitors, i.e., we obtain $S^{\prime}_k\in \Pi$ out of $S_k$ such that
$$\textbf{L}(S^{\prime}_k, \Omega) < \textbf{L}(S_k,\Omega) - \rho,$$
for some uniform $\rho>0$. Since $S$ is critical for $\textbf{L}(\cdot,\Omega)$, we have $S^{\prime}_k \in \Pi$ and $\textbf{L}(S^{\prime}_k,\Omega) < \textbf{L}(\Pi,\Omega)$, for large $k \in \mathbb{N}$. This will give us a contradiction.

To understand the deformation from $S_n$ to $S^{\prime}_n$, observe that we have two types of critical varifolds $V \in \mathcal{C}(S,\Omega)$: either
\begin{enumerate}
\item[(1)] $V$ is stationary in $M$;
\item[(2)] or $\text{spt}(V)\cap \Omega = \varnothing$.
\end{enumerate}
Let $U \subset U_1 \subset \Omega$ be open sets, with $\overline{U} \subset U_1$, $\overline{U_1}\subset \Omega$ and such that $\Omega - U$ is inside a small tubular neighborhood of $\partial \Omega$.

The deformation is done in two steps. The first is based in the fact that if $V$ is type (1), then it is not almost minimizing in small annuli centered at some $p=p(V) \in \text{spt}(V)$. Then, there exists $\varepsilon(V)>0$ such that the slices $\Sigma^k_t \in S_k$ close to $V$ can be continuously deformed, increasing the $\mathcal{H}^{n-1}$-measure by an arbitrarily small amount, to a new $\tilde{\Sigma}^k_t$, such that 
$$(\tilde{\Sigma}^k_t - \Sigma^k_t)\cup (\Sigma^k_t - \tilde{\Sigma}^k_t)  \subset a(V) \quad \text{and} \quad \mathcal{H}^{n-1}(\tilde{\Sigma}^k_t) < \mathcal{H}^{n-1}(\Sigma^k_t) - \varepsilon(V),$$
where $a(V)$ is a small annulus centered at $p(V)$. If $\Sigma^k_t$ is not close to any type (1) critical varifold $V$, put $\tilde{\Sigma}^k_t = \Sigma^k_t$. Use $\tilde{S}_k = \{\tilde{\Sigma}^k_t\}_{t\in [0,1]}$.

The second part is based in the following idea: slices $\Sigma^k_t$ close to type (2) critical varifolds have small $\mathcal{H}^{n-1}(\Sigma^k_t\cap U_1)$. If we take $a(V)$ small enough, we can suppose that we always have either $a(V)\subset U_1$ or $a(V)\subset M-\overline{U}$. Anyway, we have
\begin{equation}
\mathcal{H}^{n-1}(\tilde{\Sigma}^k_t\cap \overline{U}) \leq \mathcal{H}^{n-1}(\Sigma^k_t\cap U_1).
\end{equation}
Therefore, $\tilde{S} = \{\tilde{S}_k\}_{k\in \mathbb{N}}$ is critical with respect to $\Omega$ and the big intersecting slices $\tilde{\Sigma}^k_t$ have small $\mathcal{H}^{n-1}$-measure inside $\overline{U}$. Each such slice is replaced by a hypersurface outside $\overline{\Omega}$ via a two-steps deformation. First we deform $\tilde{\Sigma}^k_t$ inside $\Omega$ to a hypersurface not intersecting $\overline{U}$ and then we use the maximum principle to take it out of $\overline{\Omega}$. And we have a contradiction.


\section{Notation and preliminaries}\label{section.pre}

\subsection{Definitions from geometric measure theory}\label{gmt.defi}

In this section we recall some definitions and notations from Geometric Measure Theory. A standard reference is the book of Simon \cite{simon}. Sometimes we will also follow the notation of  Pitts book \cite{pitts}. Our approach follows Section $4$ in \cite{marques-neves}.

Let $(M^n,g)$ be an orientable compact Riemannian manifold. We assume that $M$ is isometrically embedded in $\mathbb{R}^L$. We denote by $B(p,r)$ the open geodesic ball in $M$ of radius $r$ and center $p\in M$. 

We denote by ${\bf I}_k(M)$ the space of $k$-dimensional integral currents in $\mathbb{R}^L$ supported in $M$; ${\mathcal Z}_k(M)$ the space of integral currents $T \in {\bf I}_k(M)$ with  $\partial T=0$; and $\mathcal{V}_k(M)$ the closure of the space of $k$-dimensional rectifiable varifolds in $\mathbb{R}^L$ with support contained in $M$, in the weak topology.

Given $T\in {\bf I}_k(M)$, we denote by $|T|$ and $||T||$ the integral varifold and Radon measure in $M$ associated with $T$, respectively; given $V\in \mathcal{V}_k(M)$, $||V||$ denotes the Radon measure in $M$ associated with $V$.

The above spaces come with several relevant metrics. We use the standard notations $\mathcal{F}$ and $\textbf{M}$ for the flat norm and mass norm on ${\bf I}_k(M)$, respectively. The ${\bf F}$-{\it metric} on $\mathcal{V}_k(M)$ is defined in {Pitts book} \cite{pitts}. It induces the varifold weak topology on $\mathcal{V}_k(M)$. Finally, the ${\bf F}$-{\it metric} on ${\bf I}_k(M)$ is defined by $$ {\bf F}(S,T)=\mathcal{F}(S-T)+{\bf F}(|S|,|T|).$$ We use also the $k$-dimensional Hausdorff measure $\mathcal{H}^k$ of subsets of $M$.

We assume that ${\bf I}_k(M)$ and ${\mathcal Z}_k(M)$ both have the topology induced by the flat metric. When endowed with a different topology, these spaces will be denoted either by ${\bf I}_k(M;{\bf M})$ and ${\mathcal Z}_k(M;{\bf M})$, in case of the mass norm, or ${\bf I}_k(M;{\bf F})$ and ${\mathcal Z}_k(M;{\bf F})$, if we use the $\textbf{F}$-metric. If $U\subset M$ is an open set of finite perimeter, the associated current in ${\bf I}_{n}(M)$ is denoted by $[|U|]$. The space $\mathcal{V}_k(M)$ is always considered with the weak topology of varifolds.

Given a $C^1$-map $F:M\rightarrow M$, the push-forwards of $V\in \mathcal{V}_k(M)$ and $T\in {\bf I}_k(M)$ are denoted by  ${F}_{\#}(V)$ and ${F}_{\#}(T)$, respectively. Let $\mathcal{X}(M)$ denote the space of smooth vector fields of $M$ with the $C^1$-topology. The first variation $\delta:\mathcal{V}_k(M)\times \mathcal{X}(M)\rightarrow \R $ is defined as
$$\delta V(X)=\left.\frac{d}{dt}\right|_{t=0}||{F_{t}}_{\#}(V)||(M),$$
where $\{F_t\}_t$ is the flow of $X$. The first variation is continuous  with respect to the product topology of $\mathcal{V}_k(M)\times \mathcal{X}(M)$. A varifold $V$ is said to be {\it stationary in $M$} if $\delta V(X)=0$ for every $X\in \mathcal{X}(M)$.

\subsection{The maximum principle}\label{subsect_mp} In \cite{white}, White proved a maximum principle type theorem for general varifolds. In this paper we use an application of this result which we describe now.

Let $(M^n,g)$ be a closed Riemannian manifold and $\Omega \subset M$ be an open domain with smooth and strictly mean-concave boundary $\partial \Omega$. Consider the foliation of a small neighborhood of $\partial \Omega$ given by $\partial\Omega_t$, where $\Omega_t =\{x \in M: d_{\partial \Omega}(x) < t\}$, $d_{\partial \Omega}$ is the signed distance function to $\partial \Omega$ negative in $\Omega$ and let $\nu$ be the unit vector field normal to all $\partial\Omega_t$ pointing outside $\Omega$. Use $\overrightarrow{H}(p)$ to denote the mean-curvature vector of $\partial \Omega$ at $p \in \partial \Omega$. White's result and a covering argument in $\partial \Omega$ imply the following statement:

\subsubsection{\bf Corollary}\label{convex.slab}
\textit{Let $\overrightarrow{H}(p) = H(p) \nu(p)$ and suppose that $H> \eta >0$ over $\partial \Omega$. There exist $a<0<b$ and a smooth vector field $X$ on $M -  \Omega_a$ with
\begin{equation*}
X \cdot \nu > 0 \text{ on } \overline{\Omega_b} -  \Omega_a
\end{equation*}
and
\begin{equation*}
\delta V(X) \leq -\eta \int |X| d\mu_V,
\end{equation*}
for every $(n-1)$-varifold $V$ in $M -  \Omega_a$.
}

We refer to $X$ as the maximum principle vector field. Up to multiplying $X$ by a constant, we can suppose that its flow $\{\Phi(s,\cdot)\}_{s\geq 0}$ is such that
\begin{equation}\label{time.1.is.outside}
\Phi(1,M  -  \Omega_a) \subset M  -  \overline{\Omega_b}.
\end{equation}
The key property of this flow is that it is mass-decreasing for $(n-1)$-varifolds.

\subsubsection{Remark}\label{mass.decreasing.for.n-varifolds} Actually, $\{\Phi(s,\cdot)\}_{s\geq 0}$ is also mass-decreasing for $n$-varifolds.

\subsection{Cell complexes}\label{subsect_cell complexes} We begin by introducing the domains of our discrete maps. More details can be found in \cite{marques-neves} or \cite{pitts}.

\begin{itemize}
\item $I^n = [0,1]^n \subset \mathbb{R}^n$ and $I^n_0 = \partial I^n = I^n -  (0,1)^n$;
\item for each $j \in \mathbb{N}$, $I(1,j)$ denote the cell complex of $I^1$ whose $0$-cells and $1$-cells are, respectively, $[0], [3^{-j}], \ldots, [1-3^{-j}]$,$[1]$ and $[0,3^{-j}]$, $[3^{-j},2\cdot 3^{-j}], \ldots, [1-3^{-j},1]$;
\item $I(n,j) = I(1,j) \otimes \cdot\cdot\cdot \otimes I(1,j)$, $n$ times;
\item $I(n,j)_p = \{ \alpha_1 \otimes \cdot\cdot\cdot \otimes \alpha_n:\alpha_i \in I(1,j) \text{ and } \sum_{i=1}^n \text{dim}(\alpha_i) = p\}$;
\item $I_0(n,j)_p = I(n,j)_p\cap I^n_0$, are the $p$-cells in the boundary;
\item $\partial : I(n,j) \rightarrow I(n,j)$, the boundary homomorphism is defined by
\begin{equation*}
\partial (\alpha_1 \otimes \cdot\cdot\cdot \otimes \alpha_n) = \sum_{i=1}^n (-1)^{\sigma(i)} \alpha_1 \otimes \cdot\cdot\cdot \otimes \partial \alpha_i \otimes \cdot\cdot\cdot \otimes \alpha_n,
\end{equation*}
where $\sigma(i) = \sum_{j<i} \text{dim}(\alpha_i)$, $\partial [a,b] = [b] - [a]$ and $\partial [a] = 0$;
\item $\textbf{d} : I(n,j)_0 \times I(n,j)_0 \rightarrow \mathbb{N}$, is the grid distance, it is given by
\begin{equation*}
\textbf{d}(x,y) = 3^j\sum_{i=1}^n |x_i - y_i|;
\end{equation*}
\item $\textbf{n}(i,j): I(n,i)_0 \rightarrow I(n,j)_0$, the nearest vertex map satisfies
\begin{equation*}
\textbf{d}(x,\textbf{n}(i,j)(x)) = \min\{ \textbf{d}(x,y) : y \in I(n,j)_0\}.
\end{equation*}
\end{itemize}

\subsection{Discrete maps into $\mathcal{Z}_{n-1}(M)$ and generalized homotopies}\label{subsect_discrete.maps} Let $(M^{n},g)$ denote an orientable compact Riemannian manifold.

\subsubsection{\bf Definition} Given $\phi : I(n,j)_0 \rightarrow \mathcal{Z}_{n-1}(M)$, we define its \textit{fineness} as
\begin{equation*}
\textbf{f}(\phi) = \sup\left\{ \frac{\textbf{M}(\phi(x) - \phi(y))}{\textbf{d}(x,y)} : x,y \in I(n,j)_0, x\neq y\right\}. 
\end{equation*}

\subsubsection{Remark} If we check that $\textbf{M}(\phi(x)-\phi(y))< \delta$, for every $\textbf{d}(x,y)=1$, then we can conclude directly that $\textbf{f}(\phi) < \delta$.

\subsubsection{Definition} Let $\phi_i : I(1,k_i)_0 \rightarrow \mathcal{Z}_{n-1}(M)$, $i=1,2$, be given discrete maps. We say that \textit{$\phi_1$ is $1$-homotopic to $\phi_2$ in $(\mathcal{Z}_{n-1}(M;\textbf{M}),\{0\})$ with fineness $\delta$} if we can find $k \in \mathbb{N}$ and a map
\begin{equation*}
\psi : I(1,k)_0 \times I(1,k)_0 \rightarrow \mathcal{Z}_{n-1}(M)
\end{equation*}
with the following properties:
\begin{enumerate}
\item[(i)] $\textbf{f}(\psi) <\delta$;
\item[(ii)] $\psi([i-1],x) = \phi_i(\textbf{n}(k,k_i)(x))$, $i=1,2$ and $x\in I(1,k)_0$;
\item[(iii)] $\psi(\tau,[0]) = \psi(\tau,[1])=0$, for $\tau \in I(1,k)_0$.
\end{enumerate}


\section{Min-max theory for intersecting slices}\label{min-max.section}

In this section we describe the min-max theory that we use to prove our main results. The set up that we follow to develop that is similar to the original one introduced by Almgren and Pitts. The crucial difference is that we see only the slices intersecting a fixed closed subset $\overline{\Omega} \subset M$, which is a manifold with boundary. The aim with this is to produce an embedded closed minimal hypersurface intersecting the given domain.

Let $(M^{n},g)$ be an orientable closed Riemannian manifold and $\Omega\subset M$ be a connected open subset. We begin with the basic definitions following Almgren, Pitts and Marques-Neves, see \cite{alm2}, \cite{marques-neves}, \cite{marques-neves-2} and \cite{pitts}. 

\subsection{Definition} An \textit{$$(1,\textbf{M})-\text{homotopy sequence of mappings into } (\mathcal{Z}_{n-1}(M;\textbf{M}),\{0\})$$}
is a sequence of maps $\{\phi_i\}_{i\in \mathbb{N}}$
\begin{equation*}
\phi_i : I(1,k_i)_0 \rightarrow \mathcal{Z}_{n-1}(M),
\end{equation*}
such that $\phi_i$ is $1$-homotopic to $\phi_{i+1}$ in $(\mathcal{Z}_{n-1}(M;\textbf{M}),\{0\})$ with fineness $\delta_i$ and
\begin{enumerate}
\item[(i)] $\lim_{i\rightarrow \infty} \delta_i =0$;
\item[(ii)] $\sup\{\textbf{M}(\phi_i(x)) : x \in \text{dmn}(\phi_i) \text{ and } i \in \mathbb{N}\}< \infty$.
\end{enumerate}

The notion of homotopy between two $(1,\textbf{M})$-homotopy sequences of mappings into $(\mathcal{Z}_{n-1}(M;\textbf{M}),\{0\})$, is the following: 

\subsection{Definition} We say that $S^1 = \{\phi^1_i\}_{i\in \mathbb{N}}$ \textit{is homotopic with} $S^2 = \{\phi^2_i\}_{i\in \mathbb{N}}$ if $\phi^1_i$ is $1$-homotopic to $\phi^2_i$ with fineness $\delta_i$ and $\lim_{i\rightarrow \infty} \delta_i =0$.

One checks that to be "homotopic with" is an equivalence relation on the set of $(1,\textbf{M})$-homopoty sequences of mappings into $(\mathcal{Z}_{n-1}(M;\textbf{M}),\{0\})$. An equivalence class is called a $(1,\textbf{M})$-homotopy class of mappings into $(\mathcal{Z}_{n-1}(M;\textbf{M}),\{0\})$. We follow the usual notation $\pi^{\#}_1(\mathcal{Z}_{n-1}(M;\textbf{M}),\{0\})$ for the set of homotopy classes. These definitions are the same as in Pitts, \cite{pitts}.

The key difference is that our width consider only slices intersecting $\overline{\Omega}$. Given a map $\phi : I(1,k)_0 \rightarrow \mathcal{Z}_{n-1}(M)$ we defined its reduced domain by
\begin{equation}\label{reduced.domain}
\text{dmn}_{\Omega}(\phi) = \{x \in I(1,k)_0 : \text{spt}(||\phi(x)||)\cap \overline{\Omega} \neq \varnothing\}.
\end{equation}

\subsection{Definition}\label{width}
Let $\Pi \in \pi^{\#}_1(\mathcal{Z}_{n-1}(M;\textbf{M}),\{0\})$ be a homotopy class and $S =\{\phi_i\}_{i\in \mathbb{N}} \in \Pi$. We define
\begin{equation}
\textbf{L}(S,\Omega) = \limsup_{i\rightarrow \infty} \max\{\textbf{M}(\phi_i(x)) : x \in \text{dmn}_{\Omega}(\phi_i)\}.
\end{equation}
The \textit{width of $\Pi$ with respect to $\Omega$} is the minimum $\textbf{L}(S,\Omega)$ among all $S \in \Pi$, 
\begin{equation}
\textbf{L}(\Pi,\Omega) = \inf\{\textbf{L}(S,\Omega) : S \in \Pi\}. 
\end{equation}

Keeping the notation in the previous definition, we write $V \in \textbf{K}(S,\Omega)$ if $V = \lim_{j} |\phi_{i_j}(x_j)|$, for some increasing sequence $\{i_j\}_{j\in \mathbb{N}}$ and $x_j \in \text{dmn}_{\Omega}(\phi_{i_j})$. Moreover, if $\textbf{L}(S,\Omega) = \textbf{L}(\Pi,\Omega)$, we say that \textit{$S$ is critical with respect to $\Omega$}. In this case we consider the \textit{critical set of $S$ with respect to $\Omega$}, defined by
\begin{equation}
\mathcal{C}(S,\Omega) = \{V \in \textbf{K}(S,\Omega) : ||V||(M) = \textbf{L}(S,\Omega)\}.
\end{equation}
As in the classical theory, $\mathcal{C}(S,\Omega) \subset \mathcal{V}_{n-1}(M)$ is compact and non-empty, but here it is not clear whether there exists $V \in \mathcal{C}(S,\Omega)$ with $||V||(\overline{\Omega}) >0$.

In the direction of proving that there are critical varifolds intersecting $\overline{\Omega}$ we construct a deformation process in Section \ref{sect.small.mass}, to deal with discrete maps whose big slices enter $\Omega$ with very small mass. For doing this, and from now on, we introduce our main geometric assumption:
\begin{equation}
\Omega \textit{ has smooth and strictly mean-concave boundary } \partial \Omega.
\end{equation} 
This deformation process is inspired by Pitts' deformation arguments for constructing replacements, Section $3.10$ in \cite{pitts}, and Corollary \ref{convex.slab}. 

We introduce some notation to explain that result, which is precisely stated and proved in Section \ref{sect.small.mass}. Following the notation in the Subsection \ref{subsect_mp}, consider an open subset $U \subset \Omega$ such that $\overline{\Omega_a} \subset U$. Lemma \ref{lemma.A} guarantee the existence of positive constants $\eta_0$ and $\varepsilon_2$ depending on $M$, $\Omega$ and $U$, and $C_1$ depending only on $M$, with the following properties: given a discrete map $\phi : I(1,k)_0\rightarrow \mathcal{Z}_{n-1}(M)$ such that $\textbf{f}(\phi) \leq \eta_0$ and, for some $L>0$,
\begin{equation*}
\textbf{M}(\phi(x))\geq L \Rightarrow ||\phi(x)||(U) < \varepsilon_2,
\end{equation*}
then, up to a discrete homotopy of fineness $C_1 \textbf{f}(\phi)$, we can suppose that
\begin{equation*}
\max\{\textbf{M}(\phi(x)) : x \in \text{dmn}_{\Omega}(\phi)\} < L + C_1 \textbf{f}(\phi).
\end{equation*}
The proof is quite technical and lengthy, because it involves interpolation results, see Section \ref{interp.cont.to.discrete}. Despite the technical objects in that proof, the lemma has several applications in key arguments of this work.

It is important to generate $(1,\textbf{M})$-homotopy sequences of mappings into $(\mathcal{Z}_{n-1}(M;\textbf{M}),\{0\})$ out of a continuous map $\Gamma: [0,1]\rightarrow \mathcal{Z}_{n-1}(M;\textbf{F})$, with $\Gamma(0) = \Gamma(1) = 0$. Similarly to the discrete set up, use the notations
\begin{equation}
\text{dmn}_{\Omega}(\Gamma) = \{t \in [0,1] : \text{spt}(||\Gamma(t)||)\cap \overline{\Omega} \neq \varnothing\}
\end{equation}
and
\begin{equation}
L(\Gamma,\Omega) = \sup\{ \textbf{M}(\Gamma(t)) : t \in \text{dmn}_{\Omega}(\Gamma) \}.
\end{equation}

\subsection{Theorem}\label{cont.discrete}
\textit{Let $\Gamma$ be as above and suppose that it defines a non-trivial class in $\pi_1(\mathcal{Z}_{n-1}(M;\mathcal{F}),0)$. Then, there exists a non-trivial homotopy class $\Pi\in \pi_1^{\#}(\mathcal{Z}_{n-1}(M;\textbf{M}),\{0\})$, such that $\textbf{L}(\Pi,\Omega) \leq L(\Gamma,\Omega)$.
}

\subsubsection{Remark} Proving this is a nice application of Lemma \ref{lemma.A}, combined with interpolation results. Section \ref{creating.discrete.sweepouts} is devoted to this argument.

Observe that item (ii) in the definition of $(1,\textbf{M})$-homotopy sequences of mappings into $(\mathcal{Z}_{n-1}(M;\textbf{M}),\{0\})$ requires a uniform control on the masses of all slices. When we try to minimize the width in a given homotopy class, in order to construct a critical sequence with respect to $\Omega$, there is no restriction about the non-intersecting slices. Then, it is possible that item (ii) fails in the limit. Precisely, if $\Pi \in \pi_1^{\#}(\mathcal{Z}_{n-1}(M;\textbf{M}),\{0\})$ is a homotopy class, via a diagonal sequence argument through a minimizing sequence $\{S^j\}_{j\in \mathbb{N}} \subset \Pi$, we can produce $S^{\ast} = \{\phi_i^*\}_{i\in \N}$, so that $\phi_i^{\ast}$ is $1$-homotopic to $\phi_{i+1}^{\ast}$ with fineness tending to zero and
\begin{equation*}
\lim_{i\rightarrow \infty} \max\{\textbf{M}(\phi_i^*(x)) : x \in \text{dmn}_{\Omega}(\phi_i^*)\} = \textbf{L}(\Pi,\Omega).
\end{equation*}
But we are not able to guarantee that
\begin{equation*}
\sup\{\textbf{M}(\phi_i^*(x)) : x \in \text{dmn}(\phi_i^*) \text{ and } i \in \mathbb{N}\}< +\infty.
\end{equation*}
To overcome this difficulty we prove the following:

\subsection{Lemma}\label{control of non-intersecting slices}
\textit{There exists $C = C(M,\Omega) > 0$ with the following property: given a discrete map $\phi : I(1,k)_0 \rightarrow \mathcal{Z}_{n-1}(M)$ of small fineness, we can find $\tilde{\phi} : I(1,\tilde{k})_0 \rightarrow \mathcal{Z}_{n-1}(M)$ such that:
\begin{enumerate}
\item[(a)] $\tilde{\phi}$ is $1$-homotopic to $\phi$ with fineness $C \cdot \textbf{f}(\phi)$;
\item[(b)] $\phi(\text{dmn}_{\Omega}(\phi)) = \tilde{\phi}(\text{dmn}_{\Omega}(\tilde{\phi}))$;
\item[(c)] $$\max \{\textbf{M}(\tilde{\phi}(x)) : x \in \text{dmn}(\tilde{\phi})\} \leq C \cdot \left(\max \{\textbf{M}(\phi(x)) : x \in \text{dmn}_{\Omega}(\phi)\} + \textbf{f}(\phi) \right).$$
\end{enumerate}
}

To produce a true competitor out of $S^*$, for each $i$ large enough we replace $\phi_i^*$ by another discrete map, also denoted by $\phi_i^*$, using Lemma \ref{control of non-intersecting slices}. The new $S^*$ has the same intersecting slices and, as consequence of items (a) and (c), the additional property of being an element of $\Pi$. This concludes the existence of critical $S^*$ for $\textbf{L}(\Pi,\Omega)$. The proof of Lemma \ref{control of non-intersecting slices} is based on a natural deformation of boundaries of $n$-currents in $n$-dimensional manifolds with boundary. Briefly, the deformation is the image of the given $(n-1)$-boundary via the gradient flow of a Morse function with no interior local maximum. We postpone the details to Section \ref{section_natural sweepouts}.

Actually, we prove a Pull-tight type Theorem, as Theorem $4.3$ in \cite{pitts} and Proposition $8.5$ in \cite{marques-neves}. Precisely, given $\Pi \in  \pi_1^{\#}(\mathcal{Z}_{n-1}(M;{\bf M}),0)$, we obtain: 

\subsection{Proposition}\label{pull.tight}
\textit{There exists  a critical sequence $S^* \in \Pi$. For each  critical sequence $S^*$, there exists a critical sequence $S\in\Pi$ such that
\begin{itemize}
\item $\mathcal{C}(S,\Omega)\subset \mathcal{C}(S^*,\Omega)$, up to critical varifolds $\Sigma$ with $||\Sigma||(\Omega) =0$;
\item every $\Sigma\in \mathcal{C}(S,\Omega)$ is either a stationary varifold or $||\Sigma||(\Omega)=0$.
\end{itemize}
}

In this statement, critical means critical with respect to $\Omega$. The proof is postponed to Section \ref{proof.pulltight}. In the classical set up, the pull-tight gives a critical sequence for which all critical varifolds are stationary in $M$. In our case, it is enough to know that the intersecting critical varifolds are stationary.

The existence of non-trivial classes was proved by Almgren, in \cite{alm1}. In fact, he provides an isomorphism
\begin{equation}
F : \pi_1^{\#}(\mathcal{Z}_{n-1}(M;\textbf{M}),\{0\}) \rightarrow H_{n}(M).
\end{equation}
Then, we apply the Proposition $8.2$ of \cite{marques-neves-2} to guarantee that the non-trivial classes have positive width, in the sense introduced in \ref{width}.

\subsection{Lemma}\label{positive.width}
\textit{If $\Pi \in \pi_1^{\#}(\mathcal{Z}_{n-1}(M;\textbf{M}),\{0\})$ and $\Pi \neq 0$, then $\textbf{L}(\Pi,\Omega)>0$.
}

\subsubsection{Remark} Given $S = \{\phi_i\}_{i\in \mathbb{N}} \in \Pi$, the $\phi_i$'s with sufficient large $i$ can be extended to maps $\Phi_i$ continuous in the mass norm, respecting the non-intersecting property, see Theorem \ref{interpolation}. The Proposition $8.2$ in \cite{marques-neves-2} provides a lower bound on the value of $\sup\{\textbf{M}(\Phi(\theta)\llcorner B(p,r)) : \theta \in S^1\}$, for continuous maps $\Phi : S^1 \rightarrow \mathcal{Z}_{n-1}(M)$ in the flat topology. To obtain Lemma \ref{positive.width}, apply Proposition $8.2$ for those $\Phi_i$ and a small geodesic ball $\overline{B(p,r)}\subset \Omega$.

The next step in Almgren and Pitts' program is to find a critical varifold $V \in \mathcal{C}(S,\Omega)$ with the property of being almost minimizing in small annuli. This is a variational property that enables us to approximate the varifold by arbitrarily close integral cycles which are themselves almost locally area minimizing. The key characteristic of varifolds that are almost minimizing in small annuli is regularity, which is not necessarily the case for general integral stationary varifolds.

Let $M^n$ be a compact Riemannian manifold and $U$ be an open subset of $M$. Use $B(p,r)$ and $A(p,s,r) = B(p,r) - \overline{B(p,s)}$ to denote the open geodesic ball of radius $r$ and centered at $p$ and the open annulus in $M$, respectively.

\subsection{Notation}
Given $\varepsilon, \delta>0$, consider the set $\mathcal{A}(U;\varepsilon,\delta)$ of integer cycles $T\in \mathcal{Z}_{n-1}(M)$ for which the following happens: for every finite sequence $T = T_0, T_1,\ldots, T_m \in \mathcal{Z}_{n-1}(M)$, such that
\begin{eqnarray*}
\text{spt}(T_i - T) \subset U, \quad \textbf{M}(T_i,T_{i-1}) \leq \delta\quad \text{and}\quad \textbf{M}(T_i) \leq \textbf{M}(T) + \delta,
\end{eqnarray*}
it must be true that $\textbf{M}(T_m) \geq \textbf{M}(T) - \varepsilon$.

\subsection{Definition}\label{am.varifolds}
Let $V \in \mathcal{V}_{n-1}(M)$ be a rectifiable varifold in $M$. We say that \textit{$V$ is almost minimizing in $U$} if for every $\varepsilon >0$, there exists $\delta >0$ and $T \in \mathcal{A}(U;\varepsilon,\delta)$ satisfying $\textbf{F}(V,|T|) < \varepsilon$. 

\subsubsection{Remark} This definition is basically the same as in Pitts' book, the difference is that we ask $|T|$ to be $\varepsilon$-close to $V$ in the $\textbf{F}$-metric on the whole $M$, not only in $U$. This creates no problem because our definition implies Pitts' and in the next step we prove the existence of almost minimizing varifolds in the sense of \ref{am.varifolds}. This is also observed in Remark $6.4$ of \cite{xin}.

The goal of our construction is to produce min-max minimal hypersurfaces with intersecting properties. In order to obtain this we prove the following version of the existence theorem of almost minimizing varifolds:

\subsection{Theorem}\label{existence.almost.mini.varifolds}
\textit{Let $(M^n,g)$ be a closed Riemannian manifold, $n\leq 7$, and $\Pi \in \pi^{\#}_1(\mathcal{Z}_{n-1}(M;\textbf{M}),\{0\})$ be a non-trivial homotopy class. Suppose that $M$ contains an open subset $\Omega$ with smooth and strictly mean-concave boundary. There exists an integral varifold $V$ such that
\begin{enumerate}
\item[(i)] $||V||(M) = \textbf{L}(\Pi,\Omega)$;
\item[(ii)] $V$ is stationary in $M$;
\item[(iii)] $||V||(\Omega) >0$;
\item[(iv)] for each $p\in M$, there exists a positive number $r$ such that $V$ is almost minimizing in $A(p,s,r)$ for all $0 < s < r$. 
\end{enumerate} 
}

This is similar to Pitts' Theorem $4.10$ in \cite{pitts}, the difference being that we prove that the almost minimizing and stationary varifold also intersects $\Omega$. We postpone its proof to Section \ref{sect-proof.of.am}, it is a combination of Pitts' argument and Lemma \ref{lemma.A}. This is the last preliminary result to prove Theorem B.

The argument to prove Theorem B is simple now. Lemma \ref{positive.width} provides a non-trivial homotopy class, for which we can apply Theorem \ref{existence.almost.mini.varifolds} and obtain an integral critical varifold $V$ that is almost minimizing in small annuli and intersects $\Omega$. The Pitts' regularity theory, developed in Chapters $5$ to $7$ in \cite{pitts}, guarantees that $\text{spt}(||V||)$ is an embedded smooth hypersurface.


\section{Interpolation Results}\label{sec.interp.discrt.cont}

Interpolation is an important tool for passing from discrete maps of small fineness to continuous maps in the space of integral cycles and vice-versa, the fineness and continuity being with respect to two  different topologies. This type of technique appeared already in \cite{alm1}, \cite{alm2}, \cite{pitts}, \cite{marques-neves} and \cite{marques-neves-2}. In our approach, we follow mostly the Sections $13$ and $14$ of \cite{marques-neves}. We also make a remark that is important for us concerning the supports of interpolating sequences. 

We start with conditions under which a discrete map is approximated by a continuous map in the mass norm. The main result is important to prove Lemmas \ref{lemma.A} and \ref{positive.width}, and Proposition \ref{pull.tight}.

Let $(M^n,g)$ be a closed Riemannian manifold. We observe from Corollary $1.14$ in \cite{alm1} that there exists $\delta_0 >0$, depending
only on $M$, such that for every
\begin{equation*}
\psi : I(2,0)_0 \rightarrow \mathcal{Z}_{n-1}(M)
\end{equation*}
with $\textbf{f}(\psi) < \delta_0$ and $\alpha \in I(2,0)_1$ with $\partial \alpha = [b] - [a]$, we can find $Q(\alpha) \in \textbf{I}_{n}(M)$ with
\begin{equation*}
\partial Q(\alpha) = \psi([b]) - \psi([a]) \text{ and } \textbf{M}(Q(\alpha)) = \mathcal{F}(\partial Q(\alpha)).
\end{equation*}

Let $\Omega_1$ be a connected open subset of $M$, such that $\overline{\Omega_1}$ is a manifold with boundary. The first important result in this section is:

\subsection{Theorem}\label{interpolation}
\textit{
There exists $C_0>0$, depending only on $M$, such that for every map
\begin{equation*}
\psi : I(2,0)_0 \rightarrow \mathcal{Z}_{n-1}(M)
\end{equation*}
with $\textbf{f}(\psi)< \delta_0$, we can find a continuous map in the mass norm
\begin{equation*}
\Psi : I^2 \rightarrow \mathcal{Z}_{n-1}(M;\textbf{M})
\end{equation*}
such that
\begin{enumerate}
\item[(i)] $\Psi(x) = \psi(x)$, for all $x \in I(2,0)_0$;
\item[(ii)] for every $\alpha \in I(2,0)_p$, $\Psi|_{\alpha}$ depends only on the values assumed by $\psi$ on the vertices of $\alpha$;
\item[(iii)]
\begin{equation*}
\sup \{\textbf{M}(\Psi(x) - \Psi(y)) : x,y \in I^2\} \leq C_0 \sup_{\alpha \in I(2,0)_1}\{\textbf{M}(\partial Q(\alpha))\}.
\end{equation*}
\end{enumerate}
Moreover, if $\textbf{f}(\psi)< \min\{\delta_0,\mathcal{H}^n(\Omega_1)\}$ and
\begin{equation*}
\text{spt}(||\psi(0,0)||) \cup \text{spt}(||\psi(1,0)||) \subset M  -  \overline{\Omega_1},
\end{equation*}
we can choose $\Psi$ with $\text{spt}(||\Psi(t,0)||) \subset M  -  \overline{\Omega_1}$, for all $t \in [0,1]$.
}

\begin{proof}
The first part of this result is Theorem $14.1$ in Marques and Neves, \cite{marques-neves}. There the authors sketch the proof following the work of Almgren, Section $6$ of \cite{alm1}, and ideas of Pitts, Theorem $4.6$ of \cite{pitts}. To prove our second claim, we follow that sketch. They start with $\Delta$, a differentiable triangulation of $M$. Hence, if $s \in \Delta$ then the faces of $s$ also belong to $\Delta$. Given $s, s^{\prime} \in \Delta$, use the notation $s^{\prime}\subset s$ if $s^{\prime}$ is a face of $s$. Let $U(s) = \cup_{s\subset s^{\prime}} s^{\prime}$. Let $\Omega_2$ be a small connected neighborhood of $\overline{\Omega_1}$, so that 
$$\text{spt}(||\psi(0,0)||) \cup \text{spt}(||\psi(1,0)||) \cap \overline{\Omega_2} = \varnothing.$$
Up to a refinement of $\Delta$, we can suppose that
\begin{equation*}
s \in \Delta \text{ and } U(s) \cap (M -  \Omega_2) \neq \varnothing \Rightarrow U(s) \subset M -  \overline{\Omega_1}.
\end{equation*}
Note that $\text{spt}(||Q||) \cap \Omega_2 = \varnothing$, for $Q = Q([0,1]\otimes [0])$. Otherwise, since $Q \in \textbf{I}_{n}(M)$ and $\partial Q = \psi(1,0) - \psi(0,0)$ does not intersect $\Omega_2$, we would have that
\begin{equation*}
\mathcal{H}^n(\Omega_1) < \mathcal{H}^n(\Omega_2) \leq ||Q||(\Omega_2) \leq \textbf{M}(Q) \leq \textbf{f}(\psi) \leq \mathcal{H}^n(\Omega_1).
\end{equation*}
This is a contradiction and we have $\text{spt}(||Q||) \cap \Omega_2 = \varnothing$. By the construction of $\Psi$ we know that
\begin{equation*}
\text{spt}(||\Psi(t,0)||) \subset \bigcup \{U(s) : U(s) \cap \text{spt}(||Q||) \neq \varnothing\}.
\end{equation*}
But $U(s) \cap \text{spt}(Q) \neq \varnothing$ implies that $U(s) \cap (M -  \Omega_2) \neq \varnothing$. Then, the choice of the triangulation gives us $U(s) \subset M -  \overline{\Omega_1}$. This concludes the proof.
\end{proof}

As in \cite{marques-neves}, we also use the following discrete approximation result for continuous maps. Assume we have a continuous map in the flat topology $\Phi : I^m \rightarrow \mathcal{Z}_{n-1}(M)$, with the following properties:
\begin{itemize}
\item $\Phi|_{I^m_0}$ is continuous in the $\textbf{F}$-metric
\item $L(\Phi) = \sup\{\textbf{M}(\Phi(x)) : x\in I^m\}< \infty$ 
\item $\limsup_{r\rightarrow 0} \textbf{m}(\Phi,r) = 0$,
\end{itemize}
where $\textbf{m}(\Phi,r)$ is the concentration of mass of $\Phi$ in balls of radius $r$, i.e.:

\begin{equation*}{\bf m}(\Phi,r)=\sup\{||\Phi(x)||(B(p,r)):x\in I^m, p\in M\}.
\end{equation*}

\subsection{Theorem}\label{interp.cont.to.discrete}
\textit{There exist sequences of mappings
$$\phi_i : I(m,k_i)_0 \rightarrow \mathcal{Z}_{n-1}(M) \text{ and }
\psi_i : I(1,k_i)_0 \times I(m,k_i)_0 \rightarrow \mathcal{Z}_{n-1}(M),$$
with $k_i < k_{i+1}$, $\psi_i([0],\cdot) = \phi_i(\cdot)$, $\psi_i([1],\cdot) = \phi_{i+1}(\cdot)|_{I(m,k_i)_0}$, and sequences $\{\delta_i\}_{i\in \mathbb{N}}$ tending to zero and $\{l_i\}_{i\in \mathbb{N}}$ tending to infinity, such that
\begin{enumerate}
\item[(i)] for every $y \in I(m,k_i)_0$
\begin{equation*}
\textbf{M}(\phi_i(y)) \leq \sup\{\textbf{M}(\Phi(x)) : \alpha \in I(m,l_i)_m, x,y \in \alpha\} + \delta_i.
\end{equation*}
In particular,
\begin{equation*}
\max\{\textbf{M}(\phi_i(x)) : x \in I(m,k_i)_0\} \leq L(\Phi) + \delta_i;
\end{equation*}
\item[(ii)] $\textbf{f}(\psi_i) < \delta_i$;
\item[(iii)] 
\begin{equation*}
\sup\{\mathcal{F}(\psi_i(y,x) - \Phi(x)) : (y,x) \in \text{dmn}(\psi_i)\} < \delta_i;
\end{equation*}
\item[(iv)] if $x \in I_0(m,k_i)_0$ and $y\in I(1,k_i)_0$, we have
\begin{equation*}
\textbf{M}(\psi_i(y,x)) \leq \textbf{M}(\Phi(x)) + \delta_i.
\end{equation*}
\end{enumerate}
Moreover, if $\Phi|_{\{0\}\times I^{m-1}}$ is continuous in the mass topology then we can choose $\phi_i$ so that
\begin{equation*}
\phi_i(x) = \Phi(x), \text{ for all } x \in B(m,k_i)_0.
\end{equation*}
}

\subsubsection{Remark}\label{rmk.interp.cont.F-metric} In case $\Phi$ is continuous in the $F$-metric on the whole $I^m$, there is no concentration of mass. This is the content of lemma $15.2$ in \cite{marques-neves}.

In the proof of our main lemma, in Section \ref{sect.small.mass}, we apply the following consequence of Theorem \ref{interp.cont.to.discrete}. Let $U, \Omega_1 \subset M$ be open subsets, being $\overline{\Omega_1}$ a manifold with boundary, and $\rho > 0$. 

Assume we have a continuous map in the $F$-metric $\Psi : I^2 \rightarrow \mathcal{Z}_{n-1}(M)$. Suppose also that
\begin{itemize}
\item $t \in [0,1] \mapsto \Psi(0,t)$ is continuous in the mass norm;
\item $\sup\{ ||\Psi(s,t)||(\overline{U}) : s\in [0,1] \text{ and } t \in \{0,1\} \} < \rho$;
\item $\text{spt}(||\Psi(1,t)||) \subset M - \overline{\Omega_1}$, for every $t \in [0,1]$.
\end{itemize}

\subsection{Corollary}\label{cor.interp.supports}
\textit{Given $\delta > 0$, there exists $k \in \mathbb{N}$ and a map
\begin{equation*}
\Psi_1 : I(2,k)_0 \rightarrow \mathcal{Z}_{n-1}(M),
\end{equation*}
with the following properties:
\begin{enumerate}
\item[(i)] $\textbf{f}(\Psi_1) < \delta$;
\item[(ii)] $\sup\{ ||\Psi_1(\sigma,\tau)||(\overline{U}) : \sigma \in I(1,k)_0 \text{ and } \tau \in \{0,1\} \} < \rho + \delta$;
\item[(iii)] $\sup\{||\Psi_1(1,\tau))||(\overline{\Omega_1}) : \tau \in I(1,k)_0 \}< \delta$;
\item[(iv)] $\textbf{M}(\Psi_1(x)) \leq \textbf{M}(\Psi(x)) + \delta$, for every $x \in I(2,k)_0$;
\item[(v)] $\Psi_1(0,\tau) = \Psi(0,\tau)$, if $\tau \in I(1,k)_0$.
\end{enumerate}
}

Items (i) and (v) are easy consequences of Theorem \ref{interp.cont.to.discrete}. The other items hold if we choose a sufficiently close discrete approximation. The proof of \ref{cor.interp.supports} involves a simple combination of the uniform continuity of $\textbf{M}\circ \Psi$ on $I^2$, compactness arguments and the Lemma $4.1$ in \cite{marques-neves}, that we state now. 

\subsection{Lemma}\label{famous4.1}
\textit{Let $\mathcal{S} \subset \mathcal{Z}_k(M;\textbf{F})$ be a compact set. For every $\rho > 0$, there exists $\delta$ so that for every $S \in \mathcal{S}$ and $T \in \mathcal{Z}_k(M)$
\begin{equation*}
\textbf{M}(T) < \textbf{M}(S) + \delta \text{ and } \mathcal{F}(T-S)\leq \delta \Rightarrow \textbf{F}(S,T) \leq \rho.
\end{equation*}
}

The last lemma that we discuss in this section is similar to Theorem \ref{interp.cont.to.discrete} restricted to the case $m=1$. It says that the hypothesis about the no concentration of mass is not required in this case. See Lemma $3.8$ in \cite{pitts}.

\subsection{Lemma}\label{interp.flat-mass}
\textit{Suppose $L, \eta>0$, $K$ compact subset of $U$ and $T\in \mathcal{Z}_{k}(M)$. There exists $\varepsilon_0 = \varepsilon_0(L,\eta,K,U,T)>0$, such that whenever
\begin{itemize}
\item $S_1, S_2 \in \mathcal{Z}_k(M)$;
\item $\mathcal{F}(S_1-S_2) \leq \varepsilon_0$;
\item $\text{spt}(S_1-T)\cup \text{spt}(S_2-T) \subset K$;
\item $\textbf{M}(S_1)\leq L$ and $\textbf{M}(S_2)\leq L$,
\end{itemize}
there exists a finite sequence $S_1 = T_0, T_1,\ldots,T_m = S_2 \in \mathcal{Z}_k(M)$ with 
\begin{equation*}
\text{spt}(T_l - T) \subset U, \quad \textbf{M}(T_l-T_{l-1})\leq \eta \text{ and } \textbf{M}(T_l)\leq L+\eta.
\end{equation*}
}

This lemma is useful to approximate continuous maps in flat topology by discrete ones with small fineness in mass norm, but first we have to restrict the continuous map to finer and finer grids $I(1,k)_0$.

\section{Natural deformations on manifolds with boundary}\label{section_natural sweepouts}

In this section, we give the detail to the proof of Lemma \ref{control of non-intersecting slices}. But first, we state the following deformation result for boundaries of $n$-currents in compact $n$-dimensional manifolds with boundary.

\subsection{Lemma}\label{natural_sweepout}
\textit{Let $(\tilde{M}^n,g)$ be a compact Riemannian manifold with smooth boundary $\partial \tilde{M}$. There exists $C=C(\tilde{M})>0$ with the following property: given $A \in \textbf{I}_n(\tilde{M})$ such that $\textbf{M}(A) + \textbf{M}(\partial A) < \infty$, we can find a map $$\phi : [0,1] \rightarrow Z_{n-1}(\tilde{M})$$ continuous in the flat topology such that
\begin{enumerate}
\item[(i)] $\phi(0) = \partial A$ and $\phi(1) = 0$;
\item[(ii)] $\textbf{M}(\phi(t)) \leq C\cdot \textbf{M}(\partial A)$, for every $t \in [0,1]$.
\end{enumerate}
Moreover, if $\text{spt}(||A||)\subset \text{int}(\tilde{M})$, there exists a compact subset $K \subset \text{int}(\tilde{M})$, such that $\text{spt}(||\phi(t)||) \subset K$, for every $t \in [0,1]$.
}

\subsubsection{Remark} The construction makes clear that $\phi$ has no concentration of mass, i.e. $\limsup_{r\rightarrow 0} \sup\{||\Phi(x)||(B(p,r)):x\in I^m, p\in \tilde{M}\} = 0$ .

Let us present the proof of Lemma \ref{control of non-intersecting slices} now, and later in this section we come back to the argument of Lemma \ref{natural_sweepout}.

\subsection*{Proof of Lemma \ref{control of non-intersecting slices}}

Consider a discrete map $\psi : I(1,k)_0 \rightarrow \mathcal{Z}_{n-1}(M)$ such that $\textbf{f}(\psi) \leq \min\{\delta_0,\mathcal{H}^n(\Omega)\}$. Recall the choice of $\delta_0$ in Section \ref{sec.interp.discrt.cont}. We replace $\psi$ with another discrete map of small fineness. Observe, for each $v \in I(1,k)_0  -  \{1\}$, we can find the isoperimetric choice $A(x) \in \textbf{I}_{n}(M)$ of $\psi(v+3^{-k})  - \psi(v)$, i.e., $\partial A(v) = \psi(v+3^{-k})  - \psi(v) \text{ and }\textbf{M}(A(v)) = \mathcal{F}(\partial A(v))$.

Consider $x,y \in I(1,k)_0$ with the following properties:
\begin{itemize}
\item $[x,y]\cap \text{dmn}_{\Omega}(\psi) = \varnothing$;
\item $x - 3^{-k}$ and $y + 3^{-k}$ belong to $\text{dmn}_{\Omega}(\psi)$.
\end{itemize}
Observe that, for every $z \in [x,y-3^{-k}]\cap I(1,k)_0$, we have
\begin{equation}
\text{spt}(||A(z)||) \subset M  -  \Omega.
\end{equation}
In fact, observe $\psi(z)$ and $\psi(z+3^{-k})$ have zero mass in $\overline{\Omega}$ and $\psi(z+3^{-k}) = \psi(z) + \partial A(z)$, so, the possibilities are either $\text{spt}(||A(z)||)$ contains $\Omega$ or do not intersect it. The first case is not possible because $A(z)$ is an integral current with mass smaller than $\textbf{f}(\psi)\leq \mathcal{H}^n(\Omega)$. Next, we explain how $\psi$ is modified on each such $[x,y]$, which is called a maximal interval of non-intersecting slices. Observe that $\text{spt}(||\psi(x)||)\cup \text{spt}(||\psi(y)||) \subset M -  \overline{\Omega}$. Let $\tilde{M} = M  -  \Omega$ and consider the integral current
\begin{equation*}
A = \sum A(z) \in \textbf{I}_n(\tilde{M}),
\end{equation*}
where we sum over all $z \in [x,y-3^{-k}]\cap I(1,k)_0$. Note that $\partial A = \psi(y) - \psi(x)$ and that $\text{spt}(||A||)$ do not intersect $\partial \tilde{M}$. Let $\phi : [0,1] \rightarrow \mathcal{Z}_{n-1}(\tilde{M})$ be the map obtained via Lemma \ref{natural_sweepout} applied to the chosen $A$.

To produce a discrete map, we have to discretize $\phi$. Since we are dealing with a one-parameter map, we can directly apply Lemma \ref{interp.flat-mass}. Consider a compact set $K \subset M -  \overline{\Omega}=: U$, such that
\begin{equation}
\text{spt}(||\phi(t)||) \subset K, \quad \text{for all } t\in [0,1].
\end{equation}
This is related to the extra property of $\phi$ that we have in Lemma \ref{natural_sweepout}, because $\text{spt}(||A||)\subset M -  \overline{\Omega}$. Following the notation of Lemma \ref{interp.flat-mass}, choose $\varepsilon_0>0$ making that statement work with $L = C \cdot \textbf{M}(\partial A)$, $\eta \leq \textbf{f}(\psi)$, $T=0$, and the fixed $K \subset U$. The constant $C = C(M,\Omega)>0$ is the one in Lemma \ref{natural_sweepout}. Let $k_1 \in \mathbb{N}$ be large enough, so that $\phi_1:=\phi|_{I(1,k_1)_0}$ is $\varepsilon_0$-fine in the flat topology. In this case, for all $\theta \in I(1,k_1)_0$, we have:
\begin{itemize}
\item $\mathcal{F}(\phi_1(\theta + 3^{-k_1})-\phi(\theta)) \leq \varepsilon_0$;
\item $\text{spt}(\phi(\theta)) \subset K$;
\item $\textbf{M}(\phi(\theta)) \leq L = C\cdot \textbf{M}(\partial A)$.
\end{itemize}
Lemma \ref{interp.flat-mass} says that we can take $\tilde{k} \geq k_1$, so that $\phi_1$ admits extension $\tilde{\phi}$ to $I(1,\tilde{k})_0$ with fineness $\textbf{f}(\tilde{\phi}) \leq \eta \leq \textbf{f}(\psi)$, $\text{spt}(\tilde{\phi}(\theta)) \subset M  -  \overline{\Omega}$ and with uniformly controlled masses
\begin{equation*}
\textbf{M}(\tilde{\phi}(\theta)) \leq L + \eta \leq C \cdot \textbf{M}(\partial A) + \textbf{f}(\psi).
\end{equation*}  

We replace $\psi|_{[x,y]\cap I(1,k)_0}$ with a discrete map defined in a finer grid $$\tilde{\psi} : I(1,k+\tilde{k})_0\cap [x,y] \rightarrow \mathcal{Z}_{n-1}(M),$$ such that 
\begin{equation*}
\tilde{\psi}(w) = \psi(y) - \tilde{\phi}((w - x)\cdot 3^{k}),
\end{equation*}
if $w \in I(1,k+\tilde{k})_0 \cap [x,x+3^{-k}]$, and $\tilde{\psi}(w) = \psi(y)$, otherwise. This map has the same fineness as $\tilde{\phi}$, and so $\textbf{f}(\tilde{\psi})\leq \textbf{f}(\psi)$. Also, $\tilde{\psi}$ and $\psi$ agree in $x$ and $y$, 
$$\tilde{\psi}(x) = \psi(y) - \tilde{\phi}(0) = \psi(y) - \partial A = \psi(x) \text{ and } \tilde{\psi}(y) = \psi(y).$$
Moreover, no slice $\tilde{\psi}(w)$ intersects $\overline{\Omega}$ and
\begin{eqnarray}\label{mass.control.1}
\textbf{M}(\tilde{\psi}(w)) & \leq & \textbf{M}(\psi(y)) + C \cdot \textbf{M}(\partial A) + \textbf{f}(\psi)\\
\nonumber & \leq & C \cdot \left( \max\{\textbf{M}(\psi(z)) : z \in \text{dmn}_{\Omega}(\psi)\} + \textbf{f}(\psi)\right). 
\end{eqnarray}
The second line is possible, because $x - 3^{-k}$ and $y + 3^{-k}$ belong to $\text{dmn}_{\Omega}(\psi)$. Also, the constant appearing in the last inequality is bigger than the original, but still depending only on $M$ and $\Omega$. 

Observe one can write the original $\psi$ restricted to $[x,y]\cap I(1,k)_0$, similarly to the expression that defines $\tilde{\psi}$, as $\psi(w) = \psi(y) - \hat{\phi}(w/3)$, where
\begin{equation*}
\hat{\phi} : I(1,k+1)_0 \rightarrow \mathcal{Z}_{n-1}(M)
\end{equation*}
is given by
\begin{itemize}
\item $\hat{\phi}(\theta) = \partial A$, if $\theta \in [0,x/3]\cap I(1,k+1)_0$;
\item $\hat{\phi}(\theta) = \psi(y) - \psi(3\theta)$, if $\theta \in [x/3,y/3]\cap I(1,k+1)_0$;
\item $\hat{\phi}(\theta) = 0$, if $\theta \in [y/3,1]\cap I(1,k+1)_0$.
\end{itemize}
We concatenate the inverse direction of $\tilde{\phi}$ with $\hat{\phi}$ to construct
\begin{equation}
\overline{\phi}: I(1,\overline{k}+1)_0 \rightarrow \mathcal{Z}_{n-1}(M),
\end{equation}
where $\overline{k} = \max\{\tilde{k},k+1\}$, and
\[ \overline{\phi}(\theta) = \left\{ 
  \begin{array}{l l}
   \tilde{\phi} \circ \textbf{n}(\overline{k},\tilde{k})(1-3\theta) & \quad \text{if } 0\leq \theta \leq 3^{-1} \\
   \hat{\phi} \circ \textbf{n}(\overline{k},k+1)(3\theta -1) & \quad \text{if } 3^{-1} \leq \theta \leq 2\cdot3^{-1}\\
   0 & \quad \text{if } 2\cdot3^{-1} \leq \theta \leq 1.
  \end{array} \right.\]  

\subsection{Claim}\label{apply_almgren}
\textit{$\textbf{f}(\overline{\phi})\leq \textbf{f}(\psi)$ and $\overline{\phi}$ is $1$-homotopic to zero in $(\mathcal{Z}_{n-1}(M;\textbf{M}),\{0\})$ with fineness $C_1 \cdot \textbf{f}(\psi)$.
}

The constant $C_1>0$ is uniform, in the sense that it does not depend on $\psi$. The proof of this claim finishes the argument, because the discrete homotopy between $\overline{\phi}$ and the zero map tells us that the initial map $\psi$ and the $\tilde{\psi}$ we built are $1$-homotopic with the same fineness. Claim \ref{apply_almgren} is a consequence of Almgren's Isomorphism. $\quad\quad\quad\quad\quad\quad\quad\quad\quad\quad\quad\quad\quad\quad\quad\quad\quad\quad\quad\quad\qed$

Lemma \ref{natural_sweepout} is about the construction of natural deformations starting with the boundary of a fixed $n$-dimensional integral current $A$ in a compact manifold with boundary, and contracting it continuously and with controlled masses to the zero current. This generalizes the notion of cones in $\R^n$.

\subsection*{Proof of Lemma \ref{natural_sweepout}} To make the notation simpler, in this proof we use $M$ instead of $\tilde{M}$ to denote the manifold with boundary.

\subsubsection*{\textbf{Step $1$:}} Consider a Morse function $f :M \rightarrow [0,1]$ with $f^{-1}(1)=\partial M$ and no interior local maximum. Let $C(f) =  \{p_1,\ldots,p_k\} \subset M$ be the critical set of $f$, with $c_i = f(p_i)$ and $\text{index}(f,p_i) = \lambda_i \in \{0,1,\ldots,n-1\}$. Suppose, without loss of generality, $0 = c_k< \cdot\cdot\cdot < c_1 < 1$.

We adapt Morse's fundamental theorems to construct a list of homotopies $$h_i : [0,1] \times M_i \rightarrow M_i,\quad \text{ for } i = 1,\ldots,k,$$ defined on sublevel sets $M_1 = M$ and $M_i = M^{c_{i-1} - \varepsilon_{i-1}} = \{f \leq c_{i-1} - \varepsilon_{i-1}\}$, for $i = 2,\ldots, k$, with sufficiently small $\varepsilon_i>0$. Those maps are constructed in such a way we have the following properties:
\begin{itemize}
\item $h_i$ is smooth for all $i = 1,\ldots,k$;
\item $h_i(1,\cdot)$ is the identity map of $M_i$;
\item $h_i(0,M_i)$ is contained in $M_{i+1}$ with a $\lambda_i$-cell attached, if $i \leq k-1$;
\item $h_k(0,M_k) = f^{-1}(0)$.  
\end{itemize}  
This is achieved by Theorem \ref{morse_type_theorem}, in Appendix \ref{Appendix}.

\subsubsection*{\textbf{Step $2$:}}

Recall $A \in \textbf{I}_n(M)$ and $h_0(0,M)$ is contained in the union of $M_1$ with a $\lambda_1$-cell. Since $\lambda_1 < n$, the support of $A^1 = h_1(0,\cdot)_{\#}A$ is a subset of $M_2$. Indeed, it is an integral $n$-dimensional current and the $\lambda_1$-cell has zero $n$-dimensional measure. Inductively, we observe:
\begin{equation*}
A^i := h_i(0,\cdot)_{\#}A^{i-1} \text{ has support in } M_{i+1}, \text{ for } i = 2,\ldots, k, 
\end{equation*}
with the additional notation $M_{k+1} = f^{-1}(0)$. In particular, this allows us to concatenate the images of $A$ by the sequence of homotopies. Consider $$\Phi : [0,1] \rightarrow \textbf{I}_{n}(M)$$ starting with $\Phi(0) = A$ and inductively defined by 
\begin{equation*}
\Phi(t) = h_i(i-kt,\cdot)_{\#}\Phi((i-1)/k), \quad \text{ for } \left[ (i-1)/k,i/k\right] \text { and } i = 1, \ldots, k.
\end{equation*}

\subsubsection*{\textbf{Step $3$:}}

The homotopy formula, page $139$ in \cite{simon}, tells us
\begin{equation}\label{homotopy_formula}
\Phi(t) - \Phi(s) = (h_i)_{\#}\left( [|(i-ks,i-kt)|]\times \partial \Phi((i-1)/k)\right),
\end{equation}
for $(i-1)/k \leq t \leq s \leq i/k$. In fact, the boundary term vanishes,
\begin{equation*}
\partial (h_i)_{\#}\left( [|(i-ks,i-kt)|]\times \Phi((i-1)/k)\right) = 0,
\end{equation*}
because $[|(i-ks,i-kt)|]\times \Phi((i-1)/k)$ is $(n+1)$-dimensional, while $h_i$ takes values in a $n$-dimensional space. 

\subsubsection*{\textbf{Step $4$:}}

The mass of the product current $[|(i-ks,i-kt)|]\times \partial \Phi((i-1)/k)$ is the product measure of Lebesgue's measure on $(i-ks,i-kt)$ and the mass of $\partial \Phi((i-1)/k)$. Since all homotopies $h_1,\ldots,h_k$ are smooth, we conclude that $\Phi$ is continuous in the mass norm and that there exists a positive constant $C=C(M)>0$ such that
\begin{equation*}
\textbf{M}(\partial \Phi(t)) \leq C \cdot \textbf{M}(A), \quad \text{for all } t \in [0,1].
\end{equation*}

\subsection*{\textbf{Step $5$:}}

Consider $\phi$ being the boundary map applied to $\Phi$, i.e., $\phi(t) = \partial \Phi(t)$. Since $\Phi$ is continuous in the mass topology, it follows directly $\phi$ is continuous in the flat metric. Moreover, all maps $h_i(t,\cdot)$ with $t>0$ are diffeomorphisms, this implies $\phi$ is continuous in the $F$-metric, up to finite points. 

The last claim in the statement of Lemma \ref{natural_sweepout} follows from construction, because $\text{spt}(||A||) \subset \text{int}(\tilde{M})$ implies there exists a small $\varepsilon >0$ such that $\text{spt}(||A||)$ is contained in the sublevel set $K := \{f \leq 1-\varepsilon\}$. $\quad\quad\quad\quad\quad\qed$


\section{Deforming currents with small intersecting mass }\label{small.mass.section}

Let $(M^n,g)$ be a compact Riemannian manifold isometrically embedded in Euclidean space $\mathbb{R}^L$, $n\leq 7$. Consider open subsets $W\subset \subset U$ in $M$. The main result in this section is to prove that it is possible deform an integral current $T$ with small mass in $U$ to a current $T^{\ast}$ outside $W$. The deformation here is discrete, with support in $U$, arbitrarily small fineness and the masses along the deformation sequence can not increase much. In this part the mean-concavity is not required, observe the statements do not involve $\Omega$.

\subsection{Lemma}\label{small.mass}
\textit{There exists $\varepsilon_1 > 0$ with the following property: given 
\begin{equation*}
T \in \mathcal{Z}_{n-1}(M^n) \text{ with } ||T||(U)< \varepsilon_1 \text{ and } \eta >0,
\end{equation*}
it is possible to find a sequence $T = T_1, \ldots, T_q \in \mathcal{Z}_{n-1}(M^n)$ such that
\begin{enumerate}
\item[(1)] $\text{spt}(T_l - T) \subset U$;
\item[(2)] $\textbf{M}(T_l - T_{l-1}) \leq \eta$;
\item[(3)] $\textbf{M}(T_l) \leq \textbf{M}(T) + \eta$;
\item[(4)] $\text{spt}(T_q) \subset M  -  \overline{W}$.
\end{enumerate}
}

\subsubsection{Remark}\label{dependece.epsilon.1} The constant $\varepsilon_1$ depends only on $M$ and $d_g(W,M -  U)$.

The main ingredients to prove this result are: the key element is a discrete deformation process used by Pitts in the construction of replacements, see section $3.10$ in \cite{pitts}, and the monotonicity formula. 

\subsection{Pitts' deformation argument}\label{first.deformation}
\textit{Let $K\subset U$ be subsets of $M$, with $K$ compact and $U$ open. Given $T \in \mathcal{Z}_{k}(M)$ and $\eta >0$, there exists a finite $T = T_1, \ldots, T_q \in \mathcal{Z}_{k}(M)$ such that
\begin{enumerate}
\item[(i)] $\text{spt}(T_l - T) \subset U$;
\item[(ii)] $\textbf{M}(T_l-T_{l-1})\leq \eta$;
\item[(iii)] $\textbf{M}(T_l)\leq \textbf{M}(T)+\eta$;
\item[(iv)] $\textbf{M}(T_q)\leq \textbf{M}(T)$;
\item[(v)] $T_q$ is locally area-minimizing in $\text{int}(K)$.
\end{enumerate}
}

The monotonicity formula, imply that there are $C,r_0>0$, such that given $\Sigma^k \subset M$ minimal submanifold and $p\in \Sigma$, we have $\mathcal{H}^k(\Sigma\cap B(p,r)) \geq C r^k$, for all $0 < r < r_0$. As references for the monotonicity formula see \cite{c-m} or \cite{simon}. Now, we proceed to the proof of small mass deformation.

\begin{proof}[Proof of \ref{small.mass}:]
Consider a compact subset $K$ with $W \subset \subset K \subset U$, and let
\begin{equation*}
T \in \mathcal{Z}_{n-1}(M^n) \text{ with } ||T||(U)< \varepsilon_1 \text{ and } \eta >0,
\end{equation*}
be given, $\varepsilon_1>0$ to be chosen. Apply corollary \ref{first.deformation} to these $K\subset U$, $T$ and $\eta$, to obtain a finite sequence $T = T_1, \ldots, T_q \in \mathcal{Z}_{n-1}(M)$ with properties (i)-(v) of that result, this gives (1)-(3) of Lemma \ref{small.mass} automatically. Note that (i) and (iv) imply $||T_q||(U) \leq ||T||(U)\leq \varepsilon_1$. By (v) and the dimension restriction, $n\leq 7$, we can use the regularity theory for codimension-one area-minimizing currents to conclude that $T_q$ is a stable minimal hypersurface in $\text{int}(K)$. Suppose $T_q$ is not outside $W$, take $p \in W\cap \text{spt}(T_q)$ and $0 < r = 2^{-1}\min\{r_0,d_g(W,M -  K)\}$. Then
\begin{equation*}
Cr^{n-1} \leq ||T_q||(B(p,r))\leq ||T_q||(U) \leq \varepsilon_1.
\end{equation*}
So, in order to conclude the proof, we need only to choose $\varepsilon_1< Cr^{n-1}$ 
\end{proof}


\section{Deforming bad discrete sweepouts}\label{sect.small.mass}

Let $(M^n,g)$ be a closed embedded submanifold of $\mathbb{R}^L$ and $\Omega \subset M$ be an open subset with smooth and strictly mean-concave boundary $\partial \Omega$. Recall the domains $\Omega_a \subset \Omega \subset \Omega_b$ and the maximum principle vector field $X$ that we considered in Subsection \ref{subsect_mp}. Let $U \subset \Omega$ be an open subset such that $\overline{\Omega_a} \subset U$. In this section we consider discrete maps $\phi : I(1,k)_0 \rightarrow \mathcal{Z}_{n-1}(M)$ with small fineness for which slices $\phi(x)$ with mass greater than a given $L>0$ have small mass in $U$. The goal here is to extend the construction of Section \ref{small.mass.section} and deform the map $\phi$ via a $1$-homotopy with small fineness.

In order to state the precise result, consider the constants: $C_0=C_0(M)$ and $\delta_0=\delta_0(M)$ as introduced in Section \ref{sec.interp.discrt.cont}, and $\varepsilon_1(U,\Omega_a)$ as given by Lemma \ref{small.mass} for $W= \Omega_a$ and the fixed $U$. Consider also the following combinations
\begin{itemize}
\item $(3+C_0) \eta_0 = \varepsilon_2 = \min\{\varepsilon_1(U,\Omega_a),5^{-1}\delta_0, 5^{-1}\mathcal{H}^n(\Omega_a)\}$;
\item $C_1 = 3C_0+7$.
\end{itemize}
Observe that $C_1 = C_1(M)$, but $\varepsilon_2$ and $\eta_0$ depend also on $U$ and $\Omega_a$.

Assume we have a discrete map $\phi : I(1,k)_0 \rightarrow \mathcal{Z}_{n-1}(M)$ with fineness $\textbf{f}(\phi)\leq \eta_0$ and satisfying the property that, for some $L>0$,
\begin{equation}\label{very.big.imply.small.in.U}
\textbf{M}(\phi(x)) \geq L \Rightarrow ||\phi(x)||(\overline{U}) < \varepsilon_2.
\end{equation}

\subsection{Lemma}\label{lemma.A}
\textit{
There exists $\tilde{\phi} : I(1,N)_0 \rightarrow \mathcal{Z}_{n-1}(M)$ $1$-homotopic to $\phi$ with fineness $C_1\cdot\textbf{f}(\phi)$, with the following properties:
\begin{equation*}
L(\tilde{\phi}) = \max\{\textbf{M}(\tilde{\phi}(x)) : x \in \text{dmn}_{\Omega}(\tilde{\phi})\} < L + C_1 \cdot\textbf{f}(\phi)
\end{equation*}
and such that the image of $\tilde{\phi}$ coincide with $\{\phi(x) : \textbf{M}(\phi(x))<L\}$, up to $(n-1)$-currents $T \in \mathcal{Z}_{n-1}(M)$ with $||T||(\overline{U}) \leq 2\varepsilon_2$.
}

\subsubsection{Remark}\label{smaller.epsilon} If $0< \varepsilon \leq \varepsilon_2$ and we assume that $||\phi(x)||(\overline{U}) < \varepsilon$, in (\ref{very.big.imply.small.in.U}), we can still apply the lemma. Moreover, the image of $\tilde{\phi}$ will coincide with $\{\phi(x) : \textbf{M}(\phi(x))<L\}$ up to slices with the property $||\tilde{\phi}(x)||(\overline{U}) \leq \varepsilon + C_1 \cdot \eta$.
 

\subsection{Overwiew of the proof of Lemma \ref{lemma.A}} Deforming each big slice using the small mass procedure we obtain a map $\psi$ defined in a $2$-dimensional grid. The first difficulty that arises is that the obtained map has fineness of order $\varepsilon_2$ instead of $\textbf{f}(\phi)$. We correct this using the interpolation results, see Section \ref{sec.interp.discrt.cont}. The fine homotopy that we are able to construct ends with a discrete map $\tilde{\phi}$ whose intersecting slices $\tilde{\phi}(x)$ with mass exceeding $L$ by much do not get very deep in $\Omega$, i.e., $\text{spt}(\tilde{\phi}(x)\llcorner \Omega)$ is contained in a small tubulat neighborhood of $\partial \Omega$ in $M$.   

The second part is the application of the maximum principle. The idea is to continuously deform the slices supported in $M  - \Omega_a$ of the map $\tilde{\phi}$, produced in the first step, via the flow $\{\Phi(s,\cdot)\}_{s\geq 0}$ of the maximum principle vector field, see Corollary \ref{convex.slab}. Since $\Phi(1,M-\Omega_a) \subset M - \overline{\Omega_b}$, the bad slices $\tilde{\phi}(x)$ end outside $\overline{\Omega_b}$. But this deformation is continuous only with respect to the $F$-metric and we need a map with small fineness in the mass norm. This problem is similar to the difficulty that arises in the classical pull-tight argument. We produce then a discrete version of the maximum principle deformation that is arbitrarily close, in the $F$-metric, to the original one. This correction creates one more complication, because the approximation can create very big intersecting slices with small mass in $\Omega_b$. To overcome this we apply the small mass procedure again.


\begin{proof}[\textbf{Proof of \ref{lemma.A}}]
Consider the set
\begin{equation*}
\mathcal{K} = \{x \in I(1,k)_0 : \textbf{M}(\phi(x)) \geq L\}.
\end{equation*}
Let $\alpha,\beta \in\mathcal{K}$. A subset $[\alpha,\beta]\cap I(1,k)_0$ is called a maximal interval on $\mathcal{K}$ if $[\alpha,\beta]\cap I(1,k)_0 \subset \mathcal{K}$ and $\alpha - 3^{-k}, \beta + 3^{-k} \notin \mathcal{K}$.

We describe the construction of the homotopy on each maximal interval.
Observe that for every $x \in [\alpha,\beta]\cap I(1,k)_0$, we have $||\phi(x)||(\overline{U}) < \varepsilon_2 \leq \varepsilon_1$.


Let $N_1 = N_1(\phi)$ be a positive integer so that, for each $x \in [\alpha,\beta]\cap I(1,k)_0$, we can apply Lemma \ref{small.mass} to $W=\Omega_a\subset\subset U$ and find sequences
$$\phi(x) = T(0,x), T(1,x), \ldots, T(3^{N_1},x) \in \mathcal{Z}_{n-1}(M),$$
with fineness at most $\eta = \textbf{f}(\phi)$, controlled supports and masses
\begin{equation}\label{supports.and.masses}
\text{spt}(T(l,x) - \phi(x)) \subset U \quad \text{and} \quad \textbf{M}(T(l,x)) \leq \textbf{M}(\phi(x)) + \eta,
\end{equation}
and ending with an integral cycle $T(3^{N_1},x))$ whose support is contained in $M  -  \overline{\Omega_a}$. Observe that $||\phi(x)||(\Omega) =0$ imply that $T(l,x)$ is constant $\phi(x)$. Then, we perform the first step of the deformation. Consider
$$\psi : I(1,N_1)_0 \times ([\alpha,\beta]\cap I(1,k)_0) \rightarrow \mathcal{Z}_{n-1}(M)$$
defined by $\psi(l,x) = T(l\cdot 3^{N_1},x)$, for every $ x \in [\alpha,\beta]\cap I(1,k)_0$ and $l \in I(1,N_1)_0$. It follows directly from the construction that the map $\psi$ satisfies:
\begin{equation}\label{mass.side.slices}
\sup \{\textbf{M}(\psi(l,x)) : l\in I(1,N_1)_0 \text{ and } x \in \{\alpha, \beta\} \} \leq L + 2\eta,
\end{equation}
and 
\begin{equation}\label{support.outside.Omega_a}
\text{spt}(||\psi(1,x)||) \subset M  -  \overline{\Omega_a}, \quad \text{for every } x \in [\alpha,\beta]\cap I(1,k)_0.
\end{equation}
It also an easy fact that
\begin{equation}
\textbf{M}(\psi(3^{-N_1},x) - \psi(3^{-N_1},x+3^{-k}))\leq 3\eta,
\end{equation}
for $x \in [\alpha,\beta - 3^{-k}]\cap I(1,k)_0$. Observe that follows from (\ref{supports.and.masses}) that
\begin{equation}\label{mass.in.U.psi}
\sup \{ ||\psi(l,x)||(\overline{U}) : l \in I(1,N_1)_0\} \leq ||\phi(x)||(\overline{U}) + \eta \leq \varepsilon_2 + \eta,
\end{equation}
for every $x$, and $||\psi(1,x)||(\overline{U}) \leq ||\phi(x)||(\overline{U}) \leq \varepsilon_2$. Moreover, we check that
\begin{equation}\label{estimate.fineness}
\textbf{f}(\psi) \leq 3\eta + 4\varepsilon_2 \leq 5\varepsilon_2 \leq \min\{\delta_0, \mathcal{H}^{n}(\Omega_a)\}.
\end{equation}
We have to prove that $\textbf{M}(\psi(l,x)-\psi(l,x^{\prime}))\leq 5\varepsilon_2$, for $|x - x^{\prime}| \leq 3^{-k}$. Indeed, rewrite that difference as 
\begin{eqnarray}\label{rewrite.1}
(\psi(l,x) - \phi(x)) + (\phi(x) - \phi(x^{\prime})) + (\phi(x^{\prime}) - \psi(l,x^{\prime})).
\end{eqnarray}
The first and third terms are similar and, if they are non-zero, their analysis follow the steps: $\textbf{M}(\psi(l,x) - \phi(x)) = ||\psi(l,x) - \phi(x)||(U)$, because $\psi(l,x)$ is constructed in such a way that $\text{spt}(\psi(l,x) - \phi(x)) \subset U$ and
\begin{equation*}
||\psi(l,x) - \phi(x)||(U) \leq  || \psi(l,x)||(U) + ||\phi(x)||(U) \leq \eta + 2  ||\phi(x)||(U) < \eta + 2\varepsilon_2.
\end{equation*}
In the above estimate we used only (\ref{mass.in.U.psi}) and that $||\phi(x)||(U)< \varepsilon_2$. The mass of the second term of (\ref{rewrite.1}) is at most $\eta$, then
\begin{equation*}
\textbf{M} ( \psi(l,x) - \psi(l,x^{\prime}) ) \leq 3\eta + 4\varepsilon_2 \leq 5\varepsilon_2.
\end{equation*}
The last inequality follows from the special choices of $\eta \leq \eta_0$ and $3\eta_0\leq \varepsilon_2$.


In order to get the desired fineness on the final 1-homotopy, we have to interpolate. For each $l \in I(1,N_1)_0  - \{0, 1\}$ and $x \in [\alpha,\beta-3^{-k}]\cap I(1,k)_0$, we apply Theorem \ref{interpolation} to the restriction of $\psi$ to the four corner vertices of $[l,l+3^{-N_1}]\times[x,x+3^{-k}]$. This is allowed because of the expression (\ref{estimate.fineness}). The result of this step is a continuous map in the mass norm
\begin{equation}\label{homotopy.1}
\Psi : [3^{-N_1},1]\times [\alpha,\beta] \rightarrow \mathcal{Z}_{n-1}(M),
\end{equation}
that extends $\psi$. Moreover, in the $1$-cells of the form $[l,l+3^{-N_1}]\times\{x\}$, the interpolating elements $\Psi(s,x)$ differ from $\psi(l,x)$ or $\psi(l+3^{-N_1},x)$ in the mass norm at most by a factor of $C_0\textbf{M}(\psi(l,x) - \psi(l+3^{-N_1},x))\leq C_0\eta$.  This remark imply
\begin{equation}\label{estimate.mass.of.Psi(s,x)}
\sup \{\textbf{M}(\Psi(s,x)) : s \in [3^{-N_1},1] \text{ and } x \in \{\alpha, \beta\} \} \leq L + (2+C_0)\eta,
\end{equation}
and, together with expression (\ref{mass.in.U.psi}),
\begin{equation}\label{mass.in.U.Psi}
\sup\{||\Psi(s,x)||(\overline{U}) : s \in [3^{-N_1},1] \text{ and } x \in \{\alpha, \beta\}\} \leq (\varepsilon_2 + \eta) + C_0\eta.
\end{equation}
Since $\Psi(1,\alpha) = \psi(1,\alpha)$, we have better estimates $\textbf{M}(\Psi(1,\alpha))\leq L +2 \eta$ and $||\Psi(1,\alpha)||(\overline{U})\leq \varepsilon_2$, at time $s=1$. The same holds for $\beta$. A similar assertion holds for the $1$-cells $\{l\}\times[x,x+3^{-k}]$ and gives us
\begin{equation}\label{3C_0.factor}
\textbf{M}(\Psi(3^{-N_1},t) - \psi(3^{-N_1},x)) \leq 3 C_0 \eta,
\end{equation}
for $x \in [\alpha,\beta]\cap I(1,k)_0$ and $t \in [\alpha,\beta]$ with $|t-x|\leq 3^{-k}$. Property (\ref{3C_0.factor}) implies that restrictions of $\Psi(3^{-N_1},\cdot)$ to any $[\alpha,\beta] \cap I(1,N)_0$ are $1$-homotopic to the original $\phi$ with fineness at most $(3C_0+1)\textbf{f}(\phi)$. Finally, since $\textbf{f}(\psi)\leq \mathcal{H}^n(\Omega_a)$ and $\text{spt}(||\psi(1,x)||) \subset M - \overline{\Omega_a}$, we can suppose that
\begin{equation}\label{support.Psi(1,t)}
\text{spt}(||\Psi(1,t)||) \subset M  -  \overline{\Omega_a}, \quad \text{for every } t \in [\alpha,\beta].
\end{equation}


Recall that $\{\Phi(s,\cdot)\}_s$ is the flow of the maximum principle vector field $X$ and consider the map $$(s,t) \in [1,2]\times [\alpha,\beta] \mapsto \Phi(s-1,\cdot)_{\#}(\Psi(1,t)) =: \Psi(s,t).$$
This map is continuous in the $F$-metric, because each map $\Phi(s-1,\cdot)$ is diffeomorphism. Since $\Phi(s,\cdot)$ is a mass-decreasing flow, we have 
\begin{equation}\label{mass.is non.increasing}
\textbf{M}(\Psi(s,t)) \leq \textbf{M}(\Psi(1,t)) \leq \max\{\textbf{M}(\Psi(1,t)) : t\in [\alpha,\beta]\} < \infty,
\end{equation}
for every $t\in [\alpha,\beta]$ and $s \in [1,2]$. In particular, the estimates (\ref{estimate.mass.of.Psi(s,x)}) and (\ref{mass.in.U.Psi}) hold on the whole $[3^{-N_1},2]$.
Moreover, because of (\ref{time.1.is.outside}), we have
\begin{equation}\label{time.2.is.outside}
\text{spt}(||\Psi(2,t)||) \subset M  -  \overline{\Omega_b}, \quad \text{for every } t\in [\alpha,\beta].
\end{equation}


The final homotopy must be a fine discrete map defined on a $2$-dimensional grid. In order to attain this we interpolate once more, but now via Corollary \ref{cor.interp.supports}. Let $0< \varepsilon_3 \leq \varepsilon_1(\Omega_b,\Omega)$ be chosen in such a way that we can apply Lemma \ref{small.mass} with the sets $\Omega$ and $\Omega_b$, and with $\varepsilon_3\leq \eta$. Apply Corollary \ref{cor.interp.supports} for a sufficiently small $\delta$, to obtain a number $N_2 = N_2(\phi) \geq N_1+k$ and a discrete map 
\begin{equation*}
\Psi_1 : I(1,N_2)_0 \times ([\alpha,\beta] \cap I(1,N_2)_0) \rightarrow \mathcal{Z}_{n-1}(M),
\end{equation*}
such that
\begin{enumerate}
\item[(i)] $\textbf{f}(\Psi_1) < \eta$;
\item[(ii)]
\begin{equation*}
\sup\{ ||\Psi_1(\sigma,\tau)||(\overline{U}) : \sigma \in I(1,N_2)_0 \text{ and } \tau \in \{\alpha,\beta\} \} < (\varepsilon_2 + (1+C_0)\eta) + \eta;
\end{equation*}
\item[(iii)]
\begin{equation*}
\sup\{||\Psi_1(1,\tau))||(\overline{\Omega_b}) : \tau \in [\alpha,\beta]\cap I(1,N_2)_0 \}< \varepsilon_3;
\end{equation*}
\item[(iv)]
\begin{equation*}
\sup\{\textbf{M}(\Psi_1(\sigma,x)) : \sigma \in I(1,N_2)_0 \text{ and } x \in \{\alpha, \beta\}\} \leq  L + (3 + C_0)\eta;
\end{equation*} 
\item[(v)] $\Psi_1(0,\tau) = \Psi(3^{-N_1},\tau)$, if $\tau \in [\alpha,\beta]\cap I(1,N_2)_0$.
\end{enumerate}


The next step is a second application of the Lemma \ref{small.mass}, now for the slices $\Psi_1(1,\tau)$. The choice of $\varepsilon_3$ guarantee that there exists $N_3 = N_3(\phi) \in \mathbb{N}$ and an extension of $\Psi_1$ to the discrete domain
\begin{equation*}
\text{dmn}(\Psi_1) = (I(1,N_2)_0 \cup (I(1,N_3)_0+\{1\})) \times ([\alpha,\beta] \cap I(1,N_2)_0),
\end{equation*}
where $I(1,N_3)_0+\{1\} = \{\lambda + 1 : \lambda \in I(1,N_3)_0\}$. This map has the following properties: in the same spirit as (\ref{estimate.fineness}), we have, respectively,
\begin{equation}\label{estimate.fineness.Psi_1}
\textbf{f}(\Psi_1) \leq 3\eta + 4\varepsilon_3 \leq 7\eta.
\end{equation}
Moreover, similarly to (\ref{mass.side.slices}), (\ref{support.outside.Omega_a}) and (\ref{mass.in.U.psi}), we have, respectively,
\begin{equation}\label{mass.side.Psi_1}
\sup\{\textbf{M}(\Psi_1(\sigma,\tau)) : (\sigma,\tau) \in \text{dmn}(\Psi_1) \text{ and } \tau \in \{\alpha, \beta\}\} \leq L + (4 + C_0) \eta,
\end{equation}
\begin{equation}\label{support.outside.Omega}
\text{spt}(||\Psi_1(2,\tau)||) \subset M - \overline{\Omega}, \quad \text{for every } \tau \in ([\alpha,\beta] \cap I(1,N_2)_0) 
\end{equation}
and
\begin{equation}\label{mass.in.U.Psi_1.extension}
\sup\{ ||\Psi_1(\sigma,\tau)||(\overline{U}) : (\sigma,\tau) \in \text{dmn}(\Psi_1), \tau \in \{\alpha, \beta\} \} < \varepsilon_2 + (2+C_0)\eta.
\end{equation}


The last part is the organization of the homotopy. Take $N = N(\phi)$ sufficiently large such that it is possible to define a map
\begin{equation*}
\Psi^{\prime} : I(1,N)_0\times ([\alpha-3^{-k},\beta+3^{-k}] \cap I(1,N)_0) \rightarrow \mathcal{Z}_{n-1}(M)
\end{equation*}
in the following way:
\begin{itemize}
\item if $j= 0, 1, \ldots, 3^{N_2}$ and $\tau \in ([\alpha,\beta]\cap I(1,N)_0)$,
\begin{equation*}
\Psi^{\prime}(j\cdot 3^{-N},\tau) = \Psi_1(j\cdot 3^{-N_2},\textbf{n}(N,N_2)(\tau));
\end{equation*}
\item if $j= 0, 1, \ldots, 3^{N_3}$ and $\tau \in ([\alpha,\beta]\cap I(1,N)_0)$,
\begin{equation*}
\Psi^{\prime}(3^{-N+N_2}+ j\cdot 3^{-N},\tau) = \Psi_1(1+ j \cdot 3^{-N_3},\textbf{n}(N,N_2)(\tau)); 
\end{equation*}
\item if $0 < \lambda_2 \leq \lambda_1\leq (3^{N_2}+3^{N_3})\cdot 3^{-N}$,
\begin{equation*}
\Psi^{\prime}(\lambda_1,\alpha-\lambda_2) = \Psi^{\prime}(\lambda_1-\lambda_2,\alpha)\text{ and }\Psi^{\prime}(\lambda_1,\beta+\lambda_2) = \Psi^{\prime}(\lambda_1-\lambda_2,\beta);
\end{equation*}
\item if $0 < \lambda_1 < \lambda_2\leq (3^{N_2}+3^{N_3})\cdot 3^{-N}$,
\begin{equation*}
\Psi^{\prime}(\lambda_1,\alpha-\lambda_2) = \phi(\alpha)\text{ and }\Psi^{\prime}(\lambda_1,\beta+\lambda_2) = \phi(\beta);
\end{equation*}
\item if $0\leq \lambda_1 \leq (3^{N_2}+3^{N_3})\cdot 3^{-N} < \lambda_2 \leq 3^{-k}$
\begin{equation*}
\Psi^{\prime}(\lambda_1,\alpha-\lambda_2) = \phi(\textbf{n}(N,k)(\alpha-\lambda_2))\text{ and }\Psi^{\prime}(\lambda_1,\beta+\lambda_2) = \phi(\textbf{n}(N,k)(b+\lambda_2));
\end{equation*}
\item if $(3^{N_2}+3^{N_3})\cdot 3^{-N}\leq \lambda \leq 1$ and $\tau([\alpha-3^{-k},\beta+3^{-k}] \cap I(1,N)_0)$, put
\begin{equation*}
\Psi^{\prime}(\lambda,\tau) = \Psi^{\prime}((3^{N_2}+3^{N_3})\cdot 3^{-N},\tau).
\end{equation*}
\end{itemize}
In order to obtain a $1$-homotopy we need to take $N$ such that
\begin{equation*}
(3^{N_2}+3^{N_3})\cdot 3^{-N} < \frac{1}{2}\cdot3^{-k}.
\end{equation*}

Extend $\Psi^{\prime}$ to $I(1,N)_0\times I(1,N)_0$, using the above construction near each maximal interval on $\mathcal{K}$ and putting $\Psi^{\prime}(\lambda,\tau) = \phi(n(N,k)(\tau))$ on the complement. This map is a homotopy and has fineness $\textbf{f}(\Psi^{\prime})\leq 7 \eta$. Then, the obtained map $\tilde{\phi}(\cdot) = \Psi^{\prime}(1,\cdot)$ is $1$-homotopic to the original discrete map $\phi$ with fineness at most $\max\{7,3C_0+1\}\textbf{f}(\phi)\leq C_1\textbf{f}(\phi)$. The $(3C_0+1)$ factor comes from the first deformation step, recall the expression (\ref{3C_0.factor}).

Moreover, if $x \in ([\alpha-3^{-k},\beta+3^{-k}] \cap I(1,N)_0)$ we have three possibilities for $\tilde{\phi}(x)$: it coincide either with some $\Psi_1(\sigma,\alpha)$, or $\Psi_1(\sigma,\beta)$ or $\Psi_1(2,\tau)$. By (\ref{mass.side.Psi_1}) and (\ref{support.outside.Omega}), we conclude that  if $x \in ([\alpha-3^{-k},\beta+3^{-k}] \cap I(1,N)_0)$, then either $\textbf{M}(\tilde{\phi}(x)) \leq L + C_1\textbf{f}(\phi)$ or $x \notin \text{dmn}_{\Omega}(\tilde{\phi})$. In particular, if $x \in ([\alpha-3^{-k},\beta+3^{-k}] \cap \text{dmn}_{\Omega}(\tilde{\phi}))$, then
\begin{equation}
\textbf{M}(\tilde{\phi}(x)) \leq L + (4+C_0)\eta< L + C_1\textbf{f}(\phi).
\end{equation}
If $x \in \text{dmn}_{\Omega}(\tilde{\phi})$ and $x \notin [\alpha-3^{-k},\beta+3^{-k}]$, for any $[\alpha,\beta]$ maximal on $\mathcal{K}$, then $\tilde{\phi}(x)$ also appear in $\phi$ and $\textbf{M}(\tilde{\phi}(x))<L$. This concludes the proof.
\end{proof}


\section{Construction of discrete sweepouts}\label{creating.discrete.sweepouts}

Next, we apply Lemma \ref{lemma.A} to prove Theorem \ref{cont.discrete}.

\begin{proof}[\textbf{Proof of \ref{cont.discrete}}]
Let $\phi_i$, $\psi_i$, $\delta_i$ be given by Theorem \ref{interp.cont.to.discrete} applied to the map $\Gamma$. It follows from property (iv) of Theorem \ref{interp.cont.to.discrete} and the fact that $\Gamma(0) = \Gamma(1) = 0$ that, for all $y\in I(1,k_i)_0$ and $x\in \{0,1\}$, we have
\begin{equation}\label{mass.bound.zero.interpolation}
{\bf M}(\psi_i(y,x))\leq \delta_i.
\end{equation}
Define $\bar{\psi}_i:I(1,k_i)_0\times I(1,k_i)_0\rightarrow  \mathcal{Z}_{n-1}(M)$ by $\bar{\psi}_i(y,x)=0$ if $x\in \{0,1\}$ and $\bar{\psi}_i(y,x)=\psi_i(y,x)$ otherwise. Define also $\bar{\phi}_i(x)=\bar{\psi}_i([0],x)$ for $x \in I(1,k_i)_0$. Note that ${\bf f}(\bar{\psi_i}) < 2\delta_i$, by \eqref{mass.bound.zero.interpolation} and Theorem \ref{interp.cont.to.discrete} part (ii). 

Then, we obtain $\{\bar{\phi}_i\}_{i \in \mathbb{N}}$, that is an $(1,{\bf M})$-homotopy sequence of mappings into  $(\mathcal{Z}_{n-1}(M;{\bf M}),0)$. But, we can not control its width by $L(\Gamma,\Omega)$ yet. To simplify notation, let us keep using $\phi_i$ and $\psi_i$ instead of $\bar{\phi}_i$ and $\bar{\psi}_i$.  

Since $\Gamma$ is continuous in the $\textbf{F}$-metric, $\textbf{M}\circ \Gamma$ is uniformly continuous on the whole $I$. Combine this with item (i) of Theorem \ref{interp.cont.to.discrete} to conclude that
\begin{equation}\label{bound.mass.all.pts}
\textbf{M}(\phi_i(y)) \leq  \textbf{M}(\Gamma(y)) + \frac{1}{i} + \delta_i.
\end{equation}
From Lemma \ref{famous4.1}, properties (iii) of Theorem \ref{interp.cont.to.discrete} and (\ref{bound.mass.all.pts}), we get
\begin{equation*}
\lim_{i\to\infty}\sup\{\textbf{F}(\phi_i(x) - \Gamma(x)) : x \in \text{dmn}(\phi_i)\}=0.
\end{equation*}

Recall the domains $\Omega_t$ starting with $\Omega = \Omega_0$ defined in Subsection \ref{subsect_mp}. Also $\Omega_a$, $a<0$, fixed by Corollary \ref{convex.slab}.

Let $U \subset\subset \Omega$ be an open subset with $\overline{\Omega_a} \subset U$ and $0<\varepsilon = \varepsilon_2(U,\Omega_a)$ be given by Lemma \ref{lemma.A} with respect to the subsets $U$ and $\Omega_a$. Observe that
$$\mathcal{K} = \{\Gamma(t): t\in [0,1] \text{ and } ||\Gamma(t)||(\Omega) =0\}$$ 
is compact with respect to the $\textbf{F}$-metric. Since $\overline{U} \subset \Omega$, it is possible to find $\rho >0$ so that $\textbf{F}(\mathcal{K},T)<\rho$ implies $||T||(\overline{U})< \varepsilon_2$. This construction of $\rho$ does not involve the interpolation step, so we can suppose that, for large $i \in \mathbb{N}$, 
\begin{equation}\label{rho}
\sup\{\textbf{F}(\phi_i(x) - \Gamma(x)) : x \in \text{dmn}(\phi_i)\} < \rho.
\end{equation}
If $x \in I(1,k_i)_0$ is such that $x \in \mathcal{T}(\Gamma,\Omega)$, expression (\ref{bound.mass.all.pts}) gives us
\begin{eqnarray*}
\textbf{M}(\phi_i(x)) \leq  L(\Gamma,\Omega) + \frac{1}{i} + \delta_i.
\end{eqnarray*}
Otherwise, $\Gamma(x) \in \mathcal{K}$ and (\ref{rho}) imply $||\phi_i(x)||(\overline{U})< \varepsilon_2$. If $i\in \mathbb{N}$ is sufficiently large we apply Lemma \ref{lemma.A} to obtain $\tilde{\phi}_i$ $1$-homotopic to $\phi_i$ with fineness tending to zero and such that
\begin{equation*}
\max\{\textbf{M}(\tilde{\phi}_i(x)) : x\in \text{dmn}_{\Omega}(\tilde{\phi}_i)\} \leq L(\Gamma,\Omega) + \frac{1}{i} + \delta_i + C_1 \cdot \textbf{f}(\phi_i).
\end{equation*}
\end{proof}


\section{Pull-tight argument}\label{proof.pulltight}

The classical pull-tight argument is based on the construction of an area decreasing flow letting still the stationary varifolds. Flowing all slices of a critical sequence $S^*$, we produce a better competitor $S$, for which critical varifolds are stationary in $M$. In our setting, we use a slightly different flow, because we also let unmoved the non-intersecting varifolds.

Precisely, let $S^*=\{\phi_i^*\}_{i\in \mathbb{N}}$ be a critical sequence with respect to $\Omega$. Consider the set $A_0 \subset  \mathcal{V}_{n-1}(M)$ of varifolds with $||V||(M) \leq 2C$ and with one of the following properties: either $V$ is stationary in $M$ or $||V||(\Omega) =0$. Here, $C = \sup\{\textbf{M}(\phi_i^*(x)) : i\in \mathbb{N} \text{ and } x \in \text{dmn}(\phi_i^*)\}$.

Following the same steps as in Section $15$ of \cite{marques-neves}, we get a map
$$H:[0,1]\times(\mathcal{Z}_{n-1}(M;{\bf F})\cap \{\textbf{M}\leq 2C\})\rightarrow (\mathcal{Z}_{n-1}(M;{\bf F})\cap \{\textbf{M}\leq 2C\}),$$
whose key properties are
\begin{itemize}
\item[(i)]$H$ is continuous in the product topology;
\item[(ii)] $H(t,T)=T$ for all $0\leq t\leq 1$ if $|T|\in A_0$;
\item[(iii)] $||H(1,T)||(M)< ||T||(M)$ unless $|T|\in A_0$.
\end{itemize}

Direct application of $H$ on the slices of $S^*$, as in Section $15$ of \cite{marques-neves}, does not necessarily provide a better competitor, because this involves discrete approximations. Indeed, it is possible that the approximation near very big non-intersecting slices of $\phi_i^*$ creates bad intersecting slices. To overcome this, we make the approximation very close and then we apply Lemma \ref{lemma.A}.

Let us begin the proof.

\subsection*{Proof of \ref{pull.tight}}

We follow closely the proof given in Proposition $8.5$ of \cite{marques-neves}. The first claim about existence of critical sequences is proved in Section \ref{min-max.section}. 

We concentrate now in the pull-tight deformation of a given a critical sequence $S^*\in \Pi$. Suppose $S^*=\{\phi_i^*\}_{i\in \N}$,  and set $$C = \sup\{\textbf{M}(\phi_i^*(x)): i \in \mathbb{N} \text{ and } x\in \text{dmn}(\phi_i^*)\}< + \infty.$$
Consider the following compact subsets of $\mathcal{V}_{n-1}(M)$
\begin{align*}
A &= \{V \in \mathcal{V}_{n-1}(M): ||V||(M) \leq 2C\};\\
A_0& =\{V\in A: \text{either } V\mbox{ is stationary in } M \mbox{ or } ||V||(\Omega)=0\}.
\end{align*}

Following the same steps as in Section $15$ of \cite{marques-neves}, we get a map
$$H:[0,1]\times(\mathcal{Z}_{n-1}(M;{\bf F})\cap \{\textbf{M}\leq 2C\})\rightarrow (\mathcal{Z}_{n-1}(M;{\bf F})\cap \{\textbf{M}\leq 2C\}),$$
whose key properties are
\begin{itemize}
\item[(i)]$H$ is continuous in the product topology;
\item[(ii)] $H(t,T)=T$ for all $0\leq t\leq 1$ if $|T|\in A_0$;
\item[(iii)] $||H(1,T)||(M)< ||T||(M)$ unless $|T|\in A_0$.
\end{itemize}

We now proceed to the construction of a sequence $S=\{\phi_i\}_{i\in\N}\in \Pi$ such that $\mathcal{C}(S,\Omega)\subset A_0\cap \mathcal{C}(S^*,\Omega)$, up to critical varifolds $V$ with $||V||(\Omega) =0$. Denote the domain of $\phi_i^*$ by $I(1,k_i)_0$, and let  $\delta_i={\bf f}(\phi_i^*)$. Up to subsequence, we can suppose that, for $x\in I(1,k_i)_0$, either $\textbf{M}(\phi_i^*(x)) < \textbf{L}(\Pi,\Omega)+i^{-1}$ or $\text{spt}(||\phi_i^*(x)||)\subset M -  \overline{\Omega}$.

For sufficiently large $i$, apply Theorem \ref{interpolation} to obtain continuous maps in the mass norm
$$
\bar \Omega_i:[0,1]\rightarrow  \mathcal{Z}_{n-1}(M;{\bf M}),
$$
such that for all $x\in  I(1,k_i)_0$ and $\alpha\in I(1,k_i)_1$ we have
\begin{equation}\label{massomega} 
 \bar\Omega_i(x)=\phi_i^*(x)\quad\mbox{and}\quad\sup_{y,z\in\alpha}\{{\bf M}(\bar \Omega_i(z)-\bar \Omega_i(y))\}\leq C_0\delta_i.
\end{equation}
Moreover, for every $x \in [0,1]$, either $\textbf{M}(\overline{\Omega}_i(x)) < \textbf{L}(\Pi,\Omega)+ i^{-1} +C_0 \delta_i$ or $\text{spt}(||\overline{\Omega}_i(x)||) \subset M -  \overline{\Omega}$. 

Consider the continuous map in the ${\bf F}$-metric
$$\Omega_i:[0,1] \times [0,1]\rightarrow  \mathcal{Z}_{n-1}(M;{\bf F}),\quad  \Omega_i(t,x)=H(t,\overline{\Omega}_i(x)).$$ Observe that, for every $(t,x) \in [0,1]^2$ and large $i\in \mathbb{N}$, either $\textbf{M}(\Omega_i(t,x)) \leq \textbf{M}(\overline{\Omega}_i(x)) < \textbf{L}(\Pi,\Omega)+ i^{-1} + C_0 \delta_i$ or $\text{spt}(||\Omega_i(t,x)||) \subset M -  \overline{\Omega}$.

For each $i$, the map $\Omega_i$ above has no concentration of mass, because it is continuous in the $\textbf{F}$-metric, and has uniformly bounded masses because of property (iii) of $H$. Then, we can apply Theorem \ref{interp.cont.to.discrete} for $\Omega_i$ to obtain
$$\bar \phi_{ij}: I(1,s_{ij})_0 \times I(1,s_{ij})_0 \rightarrow   \mathcal{Z}_{n-1}(M)$$
such that 
\begin{itemize}
\item[(a)] 
${\bf f}(\bar \phi_{ij})<\frac{1}{j};$
 \item[(b)]
 $$\sup\{{\mathcal{F}}(\bar \phi_{ij}(t,x)-\Omega_i(t,x))\,:\, (t,x)\in I(2,s_{ij})_0\}\leq \frac{1}{j};
$$
\item[(c)]
$$
{\bf M}(\bar \phi_{ij}(t,x))\leq {\bf M}(\Omega_i(t,x))+\frac 1j\quad\mbox{for all $(t,x)\in I_0(2,s_{ij})_0$};
$$
\item[(d)]$\bar \phi_{ij}([0],x)=\Omega_i(0,x)=\bar\Omega_i(x)$ for all  $x\in  I(1,s_{ij})_0.$
\end{itemize}

From Lemma \ref{famous4.1}, properties (b), and (c), we get
\begin{equation*}
\lim_{j\to\infty}\sup\{{\bf F}(\bar \phi_{ij}(t,x),\Omega_i(t,x)):(t,x)\in I_0(2,s_{ij})_0\}=0.
\end{equation*}
Hence, using a diagonal sequence argument, can find $\{\bar \phi_i=\bar \phi_{ij(i)}\}$ such that 
\begin{equation}\label{phi_ij.omega_ij}
\lim_{i\to\infty}\sup\{{\bf F}(\bar \phi_{i}(t,x),\Omega_i(t,x)):(t,x)\in I_0(2,s_{ij(i)})_0\}=0.
\end{equation}

We define $\hat \phi_i: I(1,s_{ij(i)})_0\times I(1,s_{ij(i)})_0\rightarrow   \mathcal{Z}_{n-1}(M)$ to be equal to zero on $I(1,s_{ij(i)})_0 \times \{0,1\}$, and equal to $\bar \phi_i$ otherwise.

Recall the domains $\Omega_t$ starting with $\Omega = \Omega_0$ defined in Subsection \ref{subsect_mp}, as well as the fixed $\Omega_a$, $a<0$. We hope this notation do not cause confusion with the maps $\Omega_i$. Here we use sets $\Omega_t$ with negative $t$ close to zero only.

Let $U_i = \Omega_{-i^{-1}}$, $\overline{\Omega_a} \subset U_i$ be such that $U_i \subset U_{i+1}$ and $U_i \subset\subset \Omega$. Choose $0< \varepsilon_i < \varepsilon_2(U_i,\Omega_a)$, where $\varepsilon_2(U_i,\Omega_a)$ is the required to apply Lemma \ref{lemma.A} with respect to the sets $U_i$ and $\Omega_a$, and such that $\lim_{i\rightarrow \infty} \varepsilon_i =0$. The sets
\begin{equation*}
{K}_i = \{\Omega_i(1,x): x\in [0,1] \text{ and } ||\Omega_i(1,x)||(\Omega) =0\}
\end{equation*} 
are compact with respect to the $\textbf{F}$-metric. Then, for each $i\in \mathbb{N}$, it is possible to find $\rho_i >0$, so that $\textbf{F}(\mathcal{K}_i,T)<\rho_i$ implies $||T||(\overline{U_i})< \varepsilon_i$. 

Since this construction of $\rho_i$ does not involve the last interpolation step we can suppose the choice of $j(i)$ has the additional properties $j(i)\geq i$ and
\begin{equation}\label{rho_i}
\sup\{{\bf F}(\hat{\phi}_{i}(t,x),\Omega_i(t,x)):(t,x)\in I_0(2,s_{ij(i)})_0\} < \rho_i.
\end{equation}
If $x \in I(1,s_{ij})_0$ is such that $\textbf{M}(\Omega_i(1,x)) < \textbf{L}(\Pi,\Omega)+ i^{-1} + C_0 \delta_i$, then item (c) above and the choice of $j(i)$ give us
\begin{eqnarray*}
\textbf{M}(\hat{\phi}_i(1,x)) \leq \textbf{M}(\Omega_i(1,x))+\frac{1}{j(i)} <  \textbf{L}(\Pi,\Omega)+\frac{2}{i}+C_0 \delta_i.
\end{eqnarray*}
Otherwise, we have $\text{spt}(||\Omega_i(1,x)||) \subset M -  \overline{\Omega}$ and $\Omega_i(1,x) \in {K}_i$. By (\ref{rho_i}), $||\hat{\phi}_i(1,x)||(\overline{U_i})< \varepsilon_i$. If $i\in \mathbb{N}$ is sufficiently large we apply Lemma \ref{lemma.A} to obtain $\phi_i$ $1$-homotopic to $\hat{\phi}_i(1,\cdot)$ with fineness tending to zero and such that
\begin{equation}\label{small.mass.in.pull-tight}
\max\{\textbf{M}(\phi_i(x)): x \in \text{dmn}_{\Omega}(\phi_i)\} \leq \textbf{L}(\Pi,\Omega)+\frac{2}{i}+C_0 \delta_i + C_1  \textbf{f}(\hat{\phi}_i(1,\cdot)).
\end{equation}
This works only for large $i$, because there exists a $\eta_0(U_i, \Omega_a)$, such that $\textbf{f}(\hat{\phi}_i(1,\cdot)) \leq \eta_0(U_i, \Omega_a)$ is required to use Lemma \ref{lemma.A}. Indeed, by remark \ref{dependece.epsilon.1}, the choice of $\eta_0$ and $U_i \subset U_{i+1}$, we conclude that $\eta_0(U_i,\Omega_a) \leq \eta_0(U_{i+1},\Omega_a)$, while $\textbf{f}(\hat{\phi}_i(1,\cdot))$ decreases to zero as $i$ tends to infinity. Then, if $i$ is sufficiently large we have $\textbf{f}(\hat{\phi}_i(1,\cdot)) \leq \eta_0(U_1,\Omega_a) \leq \eta_0(U_i,\Omega_a)$.

Since  $\textbf{f}(\hat \phi_i)$ tends to zero,  we obtain that $\phi_i$ is $1$-homotopic to $\hat \phi_i([0],\cdot)$ in $(\mathcal{Z}_{n-1}(M;{\bf M}),0)$ with fineness tending to zero. On the other hand, it follows from  \eqref{massomega} and property (d) that $\hat \phi_i([0],\cdot)$ is $1$-homotopic to $\phi^*_i$ in $(\mathcal{Z}_{n-1}(M;{\bf M}),0)$ with fineness tending to zero. Hence,  $S=\{\phi_i\}_{i\in \N}\in \Pi$. Furthermore, it follows from (\ref{small.mass.in.pull-tight}) that $S$ is critical with respect to $\Omega$.

We are left to show that $\mathcal{C}(S,\Omega)\subset A_0\cap \mathcal{C}(S^*,\Omega)$, up to varifolds with zero mass in $\Omega$. First we compare $\mathcal{C}(S,\Omega)$ with the critical set of $\{\hat \phi_i(1,\cdot)\}_i$. By Lemma \ref{lemma.A} and remark \ref{smaller.epsilon}, each $\phi_i(x)$ either coincide with some slice of $\hat \phi_i(1,\cdot)$ or $||\phi_i(x)||(U_i) < \varepsilon_i+C_1  \textbf{f}(\hat{\phi}_i(1,\cdot))$. Then, $\mathcal{C}(S,\Omega)\subset \mathcal{C}(\{\hat \phi_i(1,\cdot)\}_i,\Omega)$, up to critical varifolds with zero mass in $\Omega$. To conclude the proof, we claim that $\mathcal{C}(\{\hat \phi_i(1,\cdot)\}_i,\Omega) \subset A_0\cap \mathcal{C}(S^*,\Omega)$. Here we omit the proof of this fact, because it is exactly the same as in the end of the Section $15$ of \cite{marques-neves}.

\section{Existence of intersecting almost minimizing varifolds}\label{sect-proof.of.am}

In this section we prove Theorem \ref{existence.almost.mini.varifolds}.

\subsection*{Proof of \ref{existence.almost.mini.varifolds}:}

\subsubsection*{\textbf{Part 1:}} Let $\Omega_a \subset \Omega \subset \Omega_b$ be fixed as in Subsection \ref{subsect_mp}. Consider open subsets $\overline{\Omega_a} \subset U \subset U_1 \subset \Omega$ such that $\overline{U} \subset U_1$ and $\overline{U_1} \subset \Omega$. Let $S \in \Pi$ be given by Proposition \ref{pull.tight}. Write $S = \{\phi_i\}_{i\in \mathbb{N}}$ and let $I(1,k_i)_0 = \text{dmn}(\phi_i)$. Consider also $\varepsilon_2 = \varepsilon_2(U,\Omega_a)$ be given by Lemma \ref{lemma.A}.

As in the original work of Pitts, our argument is by contradiction, we homotopically deform $S$ to decrease its width with respect to $\Omega$. This will create a contradiction with the fact that $S$ is a critical sequence.

\subsubsection*{\textbf{Part 2:}} Observe that $\mathcal{C}_0(S,\Omega) := \{V \in \mathcal{C}(S,\Omega) : ||V||(\Omega) =0\}$ is compact in the weak sense of varifolds. Then, it is possible to find $\varepsilon >0$ such that
\begin{equation}\label{choice.of.epsilon}
T \in \mathcal{Z}_{n-1}(M) \text{ and } \textbf{F}(|T|,\mathcal{C}_0(S,\Omega))<\varepsilon \Rightarrow ||T||(\overline{U_1}) < 2^{-1}\varepsilon_2.
\end{equation}

\subsubsection*{\textbf{Part 3:}} Given $V \in \mathcal{C}(S,\Omega)$ with $||V||(\Omega) >0$, by contradiction assumption, there exists $p = p(V) \in \text{spt}(||V||)$ such that $V$ is not almost minimizing in small annuli centered at $p$. For $\mu = 1,2$, choose sufficiently small
\begin{equation}
a_{\mu}(V) = A(p,s_{\mu}, r_{\mu}) \text{ and } A_{\mu}(V) = A(p,\tilde{s}_{\mu}, \tilde{r}_{\mu}), 
\end{equation}
with the following properties:
\begin{itemize}
\item $\tilde{r}_1 > r_1 > s_1 > \tilde{s}_1 > 3 \tilde{r}_2 > 3r_2 > 3s_2 > 3\tilde{s}_2$;
\item if $p \in U_1$, then $A_{\mu}(V)\subset U_1$;
\item if $p \notin \overline{U}$, then $A_{\mu}(V) \cap \overline{U} = \varnothing$;
\item $V$ is not almost minimizing in $a_{\mu}(V)$.
\end{itemize}

Since $V$ is not almost minimizing in $a_{\mu}(V)$, for $\mu =1,2$, there exists $\varepsilon(V)>0$ with the following property: given
\begin{equation}
T \in \mathcal{Z}_{n-1}(M),\ \ \textbf{F}(|T|,V)<\varepsilon(V),\ \ \mu \in \{1,2\} \text{ and } \eta>0,
\end{equation}
we can find a finite sequence $T = T_0, T_1,\ldots, T_q \in \mathcal{Z}_{n-1}(M)$ such that
\begin{enumerate}
\item[(a)] $\text{spt}(T_l - T) \subset a_{\mu}(V)$;
\item[(b)] $\textbf{M}(T_l - T_{l-1}) \leq \eta$;
\item[(c)] $\textbf{M}(T_l) \leq \textbf{M}(T) + \eta$;
\item[(d)] $\textbf{M}(T_q) < \textbf{M}(T) - \varepsilon(V)$.
\end{enumerate}

The properties of those annuli $a_{\mu}(V)$ concerning the sets $U$ and $U_1$ make them slightly smaller than in Pitts' original choice. The aim with this is to obtain the property that $||T||(U_1)$ controls $||T_l||(\overline{U})$, for every $l$. Indeed, (a) implies that $\text{spt}(T_l-T)$ is always contained in either $U_1$ or $M - \overline{U}$. In the first case, we can use item (c) to prove that $||T_l||(\overline{U})\leq |T_l||(U_1)\leq ||T||(U_1) +\eta$. In the second case, it follows that $||T_l||(\overline{U}) = ||T||(\overline{U}) \leq ||T||(U_1)$.


\subsubsection*{\textbf{Part 4:}} $\mathcal{C}(S,\Omega)$ is compact. Take $V_1, V_2,\ldots, V_{\nu} \in \mathcal{C}(S,\Omega)$ such that
\begin{equation}\label{finite.cover.of.critical.set}
\mathcal{C}(S,\Omega) \subset \bigcup_{j=1}^{\nu} \{V \in \mathcal{V}_{n-1}(M) : \textbf{F}(V,V_j)< 4^{-1} \varepsilon(V_j)\},
\end{equation}
where $\varepsilon(V_j)$ is the one we chose in Part $3$ in the case of intersecting varifolds, $||V_j||(\Omega) >0$, or $\varepsilon(V_j) = \varepsilon$ as chosen in Part $2$, otherwise.

\subsubsection*{\textbf{Part 5:}} Let $\delta > 0$ and $N \in \mathbb{N}$ be such that: given
\begin{equation}
i\geq N,\quad x\in \text{dmn}_{\Omega}(\phi_i)\quad \text{and} \quad \textbf{M}(\phi_i(x))\geq \textbf{L}(S,\Omega)- 2\delta,
\end{equation}
there exists $f_1(x) \in \{1,\ldots,\nu\}$ with
\begin{equation}
\textbf{F}(|\phi_i(x)|,V_{f_1(x)})< 2^{-1}\varepsilon(V_{f_1(x)}).
\end{equation}
The existence of such numbers can be seen via a contradiction argument. Moreover, choose $\delta$ and $N$ satisfying the following two extra conditions
\begin{itemize}
\item $\delta \leq \min\{2^{-1}\varepsilon(V_j) : j = 1,\ldots \nu\}$;
\item $i\geq N$ implies $\textbf{f}(\phi_i) \leq \min\{\delta,2^{-1}\varepsilon_2, \delta_0\}$.
\end{itemize}
Where $\delta_0 = \delta_0(M)$ is defined in Section \ref{sec.interp.discrt.cont} and $\varepsilon_2 = \varepsilon_2(U,\Omega_a)$ in Part $1$.

\subsubsection*{\textbf{Part 6:}}

In this step of our proof, we use the construction done in Theorem $4.10$ of \cite{pitts}, Parts $7$ to $18$, to state a deformation result for discrete maps. Fix $i\geq N$. Let $\alpha, \beta \in \text{dmn}(\phi_i) = I(1,k_i)_0$ be such that
\begin{enumerate}
\item[(1)] $[\alpha, \beta]\cap \text{dmn}(\phi_i) \subset \text{dmn}_{\Omega}(\phi_i)$;
\item[(2)] if $x \in (\alpha, \beta]\cap \text{dmn}(\phi_i)$, then $\textbf{M}(\phi_i(x)) \geq \textbf{L}(S,\Omega) - \delta$;
\item[(3)] if $x \in [\alpha, \beta]\cap \text{dmn}(\phi_i)$, then $||V_{f_1(x)}||(\Omega)>0$.
\end{enumerate}

\subsection{Claim}\label{Pitts.7-18}
\textit{There exist a sequence $\{\delta_i\}_{i\geq N}$ tending to zero, $N(i) \geq k_i$ and
\begin{equation*}
\psi_i : I(1,N(i))_0 \times ([\alpha, \beta]\cap I(1,N(i))_0) \rightarrow \mathcal{Z}_{n-1}(M)
\end{equation*}
with the following properties:
\begin{enumerate}
\item[(i)] $\lim_{i\rightarrow \infty}\textbf{f}(\psi_i) =0$;
\item[(ii)] $\psi_i([0],x) = \phi_i(\textbf{n}(N(i),k_i)(x))$;
\item[(iii)] $\textbf{M}(\psi_i(1,\zeta)) < \textbf{M}(\psi_i(0,\zeta))- \delta + \delta_i$, for every $\zeta \in [\alpha, \beta]\cap I(1,N(i))_0$;
\item[(iv)] $\psi_i(j, \alpha) = \phi_i(\alpha)$, for every $j \in I(1,N(i))_0$; 
\item[(v)] $\{\psi_i(\lambda,x): \lambda \in I(1,N(i))_0\}$ describes the deformation obtained in Part $3$, starting with $T = \phi_i(x)$, supported in some $a_{\mu}(V_{f_1(x)})$ and fineness $\eta = \delta_i$, for every $x \in (\alpha,\beta]\cap \text{dmn}(\phi_i)$.
\end{enumerate}
}

\subsubsection{Remark} If we include $x = \alpha$ in the hypothesis (2), the construction yields a map $\psi_i$, for which (v) also holds for $x = \alpha$, instead of (iv).

Originally, Pitts wrote this argument using $27$ annuli. It is suggested by \cite{c-dl}, when dealing with one-parameter sweepouts it is enough to take only $2$.

\subsubsection*{\textbf{Part 7:}} Consider $\alpha, \beta \in \text{dmn}(\phi_i) = I(1,k_i)_0$, such that $[\alpha,\beta]$ is maximal for the property: if $x \in [\alpha+3^{-k_i},\beta - 3^{-k_i}]\cap \text{dmn}(\phi_i)$, then $x \in \text{dmn}_{\Omega}(\phi_i)$, $\textbf{M}(\phi_i(x))\geq \textbf{L}(S,\Omega) - \delta$ and $||V_{f_1(x)}||(\Omega) >0$.

Let $\{\delta_i\}_{i\in \mathbb{N}}$ and $N(i) \in \mathbb{N}$ be as in Claim \ref{Pitts.7-18}. Set $n_i = N(i)+k_i+1$ and  
\begin{equation}
L_i = \max\{\textbf{M}(\phi_i(x)): x \in\text{dmn}_{\Omega}(\phi_i)\} - \delta + \delta_i.
\end{equation}

\subsection{Claim}\label{extending.deform}
\textit{There exists a map 
\begin{equation*}
\psi_i : I(1,n_i)_0 \times ([\alpha,\beta]\cap I(1,n_i)_0) \rightarrow \mathcal{Z}_{n-1}(M)
\end{equation*}
with the following properties:
\begin{enumerate}
\item[(a)] $\lim_{i\rightarrow \infty}\textbf{f}(\psi_i) =0$;
\item[(b)] $\psi_i([0],\cdot) = \phi_i \circ \textbf{n}(n_i,k_i)$;
\item[(c)] $\psi_i(\lambda,\alpha) = \phi_i(\alpha)$ and $\psi_i(\lambda,\beta) = \phi_i(\beta)$ for every $\lambda \in I(1,n_i)_0$;
\item[(d)] $\max\{\textbf{M}(\psi_i(1,\zeta)): \zeta \in [\alpha + 3^{-k_i},\beta - 3^{-k_i}]\cap I(1,n_i)_0\} < L_i$.
\end{enumerate}
Moreover, if $\zeta \in [\alpha,\beta]\cap I(1,n_i)_0$ and $\textbf{M}(\psi_i(1,\zeta)) \geq L_i$, then
\begin{equation}
||\psi_i(1,\zeta)||(\overline{U}) \leq 2^{-1}\varepsilon_2 + \delta_i.
\end{equation} 
}

\subsubsection*{Proof of Claim \ref{extending.deform}:}

Apply Claim \ref{Pitts.7-18} on $[\alpha + 3^{-k_i},\beta - 3^{-k_i}]$ to obtain the map 
\begin{equation*}
\psi_i : I(1,N(i))_0 \times ([\alpha + 3^{-k_i},\beta - 3^{-k_i}]\cap I(1,N(i))_0) \rightarrow \mathcal{Z}_{n-1}(M).
\end{equation*}
In order to perform the extension of this $\psi_i$ to $I(1,n_i)_0 \times ([\alpha,\beta]\cap I(1,n_i)_0)$, we analyze the possibilities for $\phi_i(\alpha)$ and $\phi_i(\beta)$.

If $\alpha \in \text{dmn}_{\Omega}(\phi_i)$ and $||V_{f_1(\alpha)}||(\Omega) > 0$, we can apply Claim \ref{Pitts.7-18} directly on $[\alpha,\beta - 3^{-k_i}]$. Then, we have that the map $\psi_i$ is already defined on $I(1,N(i))_0 \times ([\alpha,\beta - 3^{-k_i}]\cap I(1,N(i))_0)$. Extend it to the desired domain simply by $\psi_i \circ \textbf{n}(n_i,N(i))$. Observe that the choice of $[\alpha,\beta]$ implies that 
\begin{equation*}
\textbf{L}(S,\Omega) - 2\delta \leq \textbf{M}(\phi_i(\alpha)) < \textbf{L}(S,\Omega) -\delta.
\end{equation*}
By item (iv) of Claim \ref{Pitts.7-18}, $\psi_i(\lambda,\alpha) = \phi_i(\alpha)$, for every $\lambda \in I(1,n_i)_0$. Item (iii) of the same statement implies that 
\begin{equation*}
\max\{\textbf{M}(\psi_i(1,\zeta)) : \zeta \in [\alpha, \beta - 3^{-k_i}]\cap I(1,n_i)_0\} < L_i.
\end{equation*}

Suppose now that $\alpha$ satisfies one of the following properties:
\begin{enumerate}
\item[(I)] either $\alpha \in \text{dmn}_{\Omega}(\phi_i)$ and $||V_{f_1(\alpha)}||(\Omega) = 0$;
\item[(II)] or $\alpha \notin \text{dmn}_{\Omega}(\phi_i)$.
\end{enumerate}
In both cases, the extension of $\psi_i$ on $I(1,n_i)_0 \times ([\alpha + 3^{-k_i}, \beta - 3^{-k_i}]\cap I(1,n_i)_0)$ is given by $\psi_i \circ \textbf{n}(n_i,N(i))$. We complete the extension in such a way that  
\begin{equation*}
\{\psi_i(\lambda,\zeta): \zeta\in [\alpha, \alpha + 3^{-k_i}]\cap I(1,n_i)_0\} \subset \{\psi_i(j,\alpha+3^{-k_i})\}_{j \in I(1,N(i))_0}\cup \{\phi_i(\alpha)\}.
\end{equation*}
The motivation for doing this is that $\phi_i(\alpha)$ and all $\psi_i(j,\alpha+3^{-k_i})$ have already small mass inside $\overline{U}$. Let us first prove this claim about the masses inside $\overline{U}$ and later we conclude the construction of the map $\psi_i$.

In case (I), we observe that $V_{f_1(\alpha)} \in \mathcal{C}_0(S,\Omega)$ and that
\begin{equation*}
\textbf{F}(|\phi_i(\alpha+3^{-k_i})|, V_{f_1(\alpha)}) \leq \textbf{f}(\phi_i) + \textbf{F}(|\phi_i(\alpha)|, V_{f_1(\alpha)}) < \varepsilon,
\end{equation*}
where the last estimate is a consequence of the choice of $N$ and $f_1$ in Part $5$. In particular, $\textbf{F}(|\phi_i(\alpha+3^{-k_i})|,\mathcal{C}_0(S,\Omega)) < \varepsilon$. Item (v) of Claim \ref{Pitts.7-18} and the comments in the end of Part $3$ imply that, for every $j \in I(1,N(i))_0$, 
\begin{equation}\label{small.mass.in.U.after.Pitts}
||\psi_i(j,\alpha + 3^{-k_i})||(\overline{U}) \leq 2^{-1}\varepsilon_2 + \delta_i. 
\end{equation}

Case (II) is simpler because $||\phi_i(\alpha)||(\overline{\Omega}) =0$ directly implies that
\begin{equation*}
||\phi_i(\alpha + 3^{-k_i})||(\overline{\Omega}) \leq \textbf{f}(\phi_i) \leq 2^{-1}\varepsilon_2.
\end{equation*}
Then, a similar analysis tells us that expression (\ref{small.mass.in.U.after.Pitts}) also holds in this case.

Finally, define $\psi_i$ on $I(1,n_i)_0 \times ((\alpha,\alpha + 3^{-k_i})\cap I(1,n_i)_0)$
by
\begin{equation*}
\psi_i(\lambda, \alpha + 3^{-k_i} -\zeta\cdot 3^{-n_i}) = \psi_i(\max\{0,\textbf{n}(n_i,N(i))(\lambda) - \zeta \cdot 3^{-N(i)}\}, \alpha + 3^{-k_i}).
\end{equation*}
Put $\psi_i(\lambda,\alpha) = \phi_i(\alpha)$, for every $\lambda \in I(1,n_i)_0$. Then, Claim \ref{extending.deform} holds.

\subsubsection*{\textbf{Part 8:}}

Let $\psi_i : I(1,n_i)_0 \times I(1,n_i)_0 \rightarrow \mathcal{Z}_{n-1}(M)$ be the map obtained in such a way that for each interval $[\alpha,\beta]$ as in Part $7$, its restriction to 
\begin{equation*}
I(1,n_i)_0 \times ([\alpha,\beta]\cap I(1,n_i)_0)
\end{equation*}
is the map of Claim \ref{extending.deform}, and $\psi_i(\lambda,\zeta) = (\phi_i \circ \textbf{n}(n_i,k_i))(\zeta)$, otherwise. 

Take $S^{\ast}= \{\phi_i^{\ast}\}_{i\in \mathbb{N}}$, where $\phi_i^{\ast} = \psi_i(1,\cdot)$ is defined on $I(1,n_i)_0$. Because the maps $\psi_i$ have fineness tending to zero, $S^{\ast} \in \Pi$.

\subsubsection*{\textbf{Part 9:}}

The sequence $S$ is critical with respect to $\Omega$, this implies that
\begin{equation*}
\lim_{i\rightarrow \infty} \max\{\textbf{M}(\phi_i(x)): x\in \text{dmn}_{\Omega}(\phi_i)\} = \textbf{L}(S,\Omega).
\end{equation*}
In particular, $\lim_{i\rightarrow \infty} L_i = \textbf{L}(S,\Omega) - \delta$. Let $i$ be sufficiently large, so that
\begin{equation*}
i\geq N, \quad \textbf{f}(\phi_i^{\ast})\leq \eta_0,\quad L_i < \textbf{L}(S,\Omega) - 2^{-1}\delta\quad \text{and} \quad \delta_i < 2^{-1}\varepsilon_2,
\end{equation*}
where $\eta_0 = \eta_0(U,\Omega_a)$ is given by Lemma \ref{lemma.A}. This choice implies that the maps $\phi_i^{\ast}$ have the following property:
\begin{equation}
\textbf{M}(\phi_i^{\ast}(\zeta)) \geq L_i \Rightarrow ||\phi_i^{\ast}(\zeta)||(\overline{U}) < \varepsilon_2.
\end{equation}
Apply Lemma \ref{lemma.A} to produce $\tilde{S} = \{\tilde{\phi}_i\}_{i \in \mathbb{N}}$ homotopic with $S^*$, such that
\begin{equation*}
\max\{\textbf{M}(\tilde{\phi}_i(x)) : x\in \text{dmn}_{\Omega}(\tilde{\phi}_i)\} < L_i + C_1 \textbf{f}(\phi_i^{\ast}),
\end{equation*}
where $C_1 = C_1(M)$ is also given by Lemma \ref{lemma.A}. In particular, $\tilde{S}\in \Pi$ and $\textbf{L}(\tilde{S},\Omega) \leq \textbf{L}(S,\Omega)-\delta = \textbf{L}(\Pi, \Omega) - \delta$. This is a contradiction.


\section{Proof of Theorem A}\label{exist.non-compact}

In this section we use the min-max theory for intersecting slices to prove the main result of our work, Theorem A as stated in the Introduction.

\subsection*{Proof of Theorem A} We divide the proof in five steps.


\subsubsection*{\textbf{Step $1$:}}

Let $p\in N$ be such that items (a) and (b) about the $\star_{k}$-condition holds for geodesic balls $B(p,R)$ centered at $p$, as defined in Section \ref{introduction}. Choose a Morse function $f : N \rightarrow [0,+\infty)$ and let $\{\Sigma_t\}_{t\geq 0}$ be the one-parameter sweepout of integral cycles induced by the level sets of $f$,
\begin{equation}
\Sigma_t := \partial \big( \{x \in N : f(x) < t\}\big) \in \mathcal{Z}_{n-1}(N).
\end{equation}
This family has the special property that $\mathcal{H}^{n-1}(\Sigma_t)$ is a continuous function. Since $\overline{\Omega} \subset N$ is compact, we conclude
\begin{equation}
L := L\left( \{\Sigma_t\}_{t\geq 0}\right) = \sup \{\mathcal{H}^{n-1}(\Sigma_t) : \text{spt}(||\Sigma_t||) \cap \overline{\Omega} \neq \varnothing\} < +\infty.
\end{equation}


\subsubsection*{\textbf{Step $2$:}}

Consider $t_0 > 0$ such that $\Sigma_{t_0}$ is a smooth regular level of $f$ and large enough to satisfy the following property:

\begin{equation*}
\Sigma^{n-1}\subset N \text{ connected minimal, } \Sigma \cap \overline{\Omega} \neq \varnothing \text{ and } \inf_{\partial\Sigma} f \geq t_0 \Rightarrow \mathcal{H}^{n-1}(\Sigma) \geq 2 L.
\end{equation*}
The existence of such $t_0$ is a consequence of the monotonicity formula for minimal hypersurfaces and the fact that $N$ satisfies the $\star_k$-condition for some $k\leq \frac{2}{n-2}$. We accomplish this argument in Section \ref{choice.slice.far}.


\subsubsection*{\textbf{Step $3$:}}

Let $(M^n,h)$ be a compact Riemannian manifold without boundary, containing an isometric copy of $\{f \leq t_0\}$ and such that $f$ extends to $M$ as a Morse function. We call this extension $f_1$. In Section \ref{change.metric.far} below we perform a construction that gives one possible $M$. Since $\overline{\Omega} \subset \{f < t_0\}$, we have a copy of $\Omega$ inside $M$, which we also denote $\Omega$. We can suppose that $f_1(M)= [0,1]$.


\subsubsection*{\textbf{Step $4$:}}

Let $\Gamma = \{\Gamma_t\}_{t \in [0,1]}$ be the sweepout of $M$ given by $\Gamma_t = f_1^{-1}(t)$. Consider the set of intersecting times
\begin{equation*}
\text{dmn}_{\Omega}(\Gamma) = \{t \in [0,1] : \text{spt}(||\Gamma_t||) \cap \overline{\Omega} \neq \varnothing\}.
\end{equation*} 
Observe that $\Gamma_t$ coincide with the slice $\Sigma_t$, for every $0 \leq t \leq t_0$, and that $t \notin \text{dmn}_{\Omega}(\Gamma)$, for $t_0 < t \leq 1$. In particular,
\begin{equation*}
L(\Gamma,\Omega) := \sup \{\mathcal{H}^{n-1}(\Gamma_t) : t \in \text{dmn}_{\Omega}(\Gamma)\}
\end{equation*}
coincides with the number $L$ defined in Step $1$. We apply now the Min-max Theory developed in Section \ref{min-max.section}, to produce a non-trivial closed embedded minimal hypersurface $\Sigma^{n-1} \subset M$ with $\mathcal{H}^{n-1}(\Sigma) \leq L$ and intersecting $\overline{\Omega}$.

Since our Min-max methods follow the discrete setting of Almgren and Pitts, we still have to construct out of $\Gamma$ a non-trivial homotopy class $\Pi \in \pi_{1}^{\#}(\mathcal{Z}_{n-1}(M;\textbf{M}),\{0\})$, such that $\textbf{L}(\Pi,\Omega)\leq L(\Gamma,\Omega) = L$. This is the content of Theorem \ref{cont.discrete}. We can apply this result because $\Gamma$ is continuous in the $\textbf{F}$-metric and non-trivial.


\subsubsection*{\textbf{Step $5$:}} 

The choice of $t_0$ in Step $2$ guarantees that any component of $\Sigma^{n-1}$ that intersects $\overline{\Omega}$, can not go outside $\{f \leq t_0\}$. Otherwise, this would imply that $2L \leq \mathcal{H}^{n-1}(\Sigma) \leq L$. In conclusion, any intersecting component of $\Sigma$ is a closed embedded minimal hypersuface in the open manifold $N$. 


\subsection{Choice of a slice far from $\Omega$}\label{choice.slice.far}

Let us prove that the $\star_k$-condition implies that any minimal hypersurface that intersects $\Omega$ and, at the same time, goes far from $\Omega$ have large Hausdorff measure $\mathcal{H}^{n-1}$. Let $\Sigma^{n-1} \subset N$ be a minimal hypersurface as in Step $2$ above. The main tool for this subsection is the following consequence of the monotonicity formula.

\begin{prop}\label{monotonicity.2}
For every $q\in B(p,R)$ and $0<s< R^{-\frac{k}{2}}$, we have
\begin{equation}
\mathcal{H}^{n-1}(\Sigma\cap B(q,s)) \geq \frac{\omega_{n-1}}{e^{(n-1)\sqrt{R^k}s}}\cdot s^{n-1},
\end{equation}
where $\omega_{n-1}$ is the volume of the unit ball in $\R^{n-1}$. 
\end{prop}

Consider $R_0\leq R_1$ and $l \in \mathbb{N}$, for which $\overline{\Omega} \subset B(p,R_1)$ and $B(p,R_1+l)\subset \{f < t_0\}$, $t_0$ to be chosen. Recall that $\Sigma$ intersects $\overline{\Omega}$ and it is not contained in the sublevel set $\{f < t_0\}$. Then, for every $i \in \{1,2\ldots,l\}$, there are points $q_{ij} \in \Sigma$, $j \in \{1,2,\ldots, \lfloor \sqrt{(R_1+i)^{k}}\rfloor\}$, whose distance in $N$ to $p$ are given by
\begin{equation*}
d(q_{ij},p) = R_1 + i - 1 + \frac{2j-1}{2\sqrt{(R_1+i)^k}}.
\end{equation*} 
Observe that $q_{ij} \in B(p,R_1+i)$ and that the balls $B_{ij} = B(q_{ij},2^{-1}(R_1+i)^{-\frac{k}{2}})$ are pairwise disjoint. Apply Proposition \ref{monotonicity.2} to conclude that
\begin{eqnarray}\label{lots.of.balls}
\mathcal{H}^{n-1}(\Sigma) & \geq & \displaystyle\sum_{i=1}^{l} \displaystyle\sum_{j} \mathcal{H}^{n-1}(\Sigma \cap B_{ij})\\
\nonumber & \geq & \frac{\omega_{n-1}}{(2\sqrt{e})^{n-1}} \displaystyle\sum_{i=1}^{l} \left(\lfloor (R_1+i)^{\frac{k}{2}}\rfloor \cdot (R_1+i)^{-\frac{k(n-1)}{2}}\right). 
\end{eqnarray}
Since $k \leq \frac{2}{n-2}$, if we keep $R_1$ fixed and let $l \in \mathbb{N}$ goes to infinity, then the right-hand side of expression (\ref{lots.of.balls}) also tends to infinity. Choose $l \in \mathbb{N}$ large, for which that is greater than $2 L$, where $L = L(\{\Sigma_t\}_{t\geq 0})$ is the number we considered in Step $1$ above. This concludes the argument, because we chose $t_0$ such that $B(p,R_1+l) \subset \{f < t_0\}$.

\subsection{Change the metric $g$ far from $\Omega$}\label{change.metric.far}

Since $t_0$ is a regular value of $f$, $f^{-1}([t_0,t_0 + 3\varepsilon])$ has no critical points for sufficiently small $\varepsilon >0$. Moreover, there exists a natural diffeomorphism $$\xi : f^{-1}([t_0,t_0 + 3\varepsilon]) \rightarrow \Sigma_{t_0} \times [0,3\varepsilon],$$ that identifies $f^{-1}(t_0)$ with $\Sigma_{t_0} \times\{0\}$. Suppose that $\Sigma_{t_0} \times [0,3\varepsilon]$ has the product metric $g|_{\Sigma_{t_0}}\times \mathcal{L}$, where $\mathcal{L}$ denotes the Lebesgue measure. Consider the pullback metric $g_1 = \xi^{\ast}(g|_{\Sigma_{t_0}}\times \mathcal{L})$ and choose a smooth bump function $\varphi : [t_0,t_0+3\varepsilon] \rightarrow [0,1]$, such that
\begin{itemize}
\item $\varphi(t) = 1$, for every $t \in [t_0,t_0 + \varepsilon]$, and
\item $\varphi(t) = 0$, for every $t \in [t_0 + 2\varepsilon,t_0 + 3\varepsilon]$.
\end{itemize}
On $f^{-1}([t_0,t_0+3\varepsilon])$, mix the original $g$ and the product metric $g_1$ using the smooth function $\varphi$, to obtain
\begin{equation}
h_1(x) = (\varphi \circ f) (x) g(x) + (1 - (\varphi \circ f) (x)) g_1(x). 
\end{equation}
This metric admits the trivial smooth extension $h_1 = g$ over $f^{-1}([0,t_0])$. Summarizing, we produced a Riemannian manifold with boundary
\begin{equation*}
M_1 = (\{f \leq t_0 + 3\varepsilon\},h_1),
\end{equation*}
with the following properties:
\begin{enumerate}
\item[(i)] $M_1$ contains an isometric copy of $\{f \leq t_0\}$ with the metric $g$;
\item[(ii)] near $\partial M_1 = f^{-1}(t_0 + 3\varepsilon)$, $h_1 = g_1$ is the product metric.
\end{enumerate}
Then, it is possible to attach two copies of $M_1$ via the identity map of $\partial M_1$. Doing this we obtain a closed manifold $M$ with a smooth Riemannian metric $h$, that coincides with $h_1$ on each half. Precisely,
\begin{equation*}
M = M_1 \cup_{\mathcal{I}} M_1, 
\end{equation*}
where $\mathcal{I}$ denotes the identity map of $\partial M_1$. The metric $h$ is smooth because of item (ii). Moreover, item (i) says that $M$ has an isometric copy of $\{f \leq t_0\}$. Finally, let us construct a Morse function $f_1$ on $M$ that coincides with $f$ in the first piece $M_1$. On the second half $M_1$, put $f_1 = 2(t_0 + 3\varepsilon) - f$.

\appendix

\section{}\label{Appendix}

Let $(M^n,g)$ be a compact Riemannian manifold with smooth non-empty boundary $\partial M$ and $f: M \rightarrow [0,1]$ be a Morse function on $M$ such that $f^{-1}(1)= \partial M$ and with no interior local maximum. The fundamental theorems in Morse theory describe how the homotopy type of the sublevel sets $M^a = \{x \in M : f(x)\leq a\}$ change as $a \in [0,1]$ varies, see Theorems $3.1$ and $3.2$ in \cite{milnor}. In this section we develop a slightly different approach to those results. It is important in the proof of Lemma \ref{natural_sweepout}.

\subsection{Theorem}\label{morse_type_theorem}
\textit{Let $0 \leq c < d \leq 1$, $\lambda \in \{0, 1,\ldots,n-1\}$ and $p \in f^{-1}(c)$ be an index $\lambda$ critical point of $f$. Suppose $\varepsilon >0$ is such that $f^{-1}([c-\varepsilon, d])$ contains no critical points other than $p$. Then, for sufficiently small $\varepsilon >0$, there exists a smooth homotopy
\begin{equation*}
h : [0,1] \times M^d \rightarrow M^d,
\end{equation*}
with the following properties:
\begin{enumerate}
\item[1.] $h(1,\cdot)$ is the identity map of $M^d$;
\item[2.] $h(0,M^d)$ is contained in $M^{c-\varepsilon}$ with a $\lambda$-cell attached.  
\end{enumerate}
}

\subsubsection{Remark} The difference of the above statement and the classical ones, is that here we are able to guarantee that the homotopy is smooth by relaxing the condition of $h(1,\cdot)$ being a retraction onto the whole $M^{c-\varepsilon}$ with the $\lambda$-cell attached.

\begin{proof}
By the classical statements $3.1$ and $3.2$ in \cite{milnor}, we know that there exists a smooth homotopy between the identity map of $M^d$ and a retraction of $M^d$ onto $M^{c-\varepsilon} \cup H$, i.e., the sublevel set with a handle $H$ attached. It is also observed that this set has smooth boundary. The final argument in the proof of Theorem $3.2$ in Milnor's book is a vertical projection of the handle in $M^{c - \varepsilon} \cup e^{\lambda}$, where $e^{\lambda}$ is a $\lambda$-cell contained in $H$. See diagram $7$, on page 19 of \cite{milnor}. This projection is not adequate for us because it is not smooth. We only adapt this step in that proof by defining a smooth projection. 

Following the notation in \cite{milnor}, let $u_1,\ldots,u_n$ be a coordinate system in a neighborhood $U$ of $p$ so that the identity
\begin{equation}
f = c - (u_1^2 + \ldots + u_{\lambda}^2) + (u_{\lambda + 1}^2 + \ldots + u_n^2),
\end{equation}
holds throughout $U$. We use $\xi = u_1^2 + \ldots + u_{\lambda}^2$ and $\eta = u_{\lambda + 1}^2 + \ldots + u_n^2$. For sufficiently small $\varepsilon>0$, the $\lambda$-cell $e^{\lambda}$ can be explicitly given by the points in $U$ with $\xi \leq \varepsilon$ and $\eta =0$. Consider $\delta>0$ so that the image of $U$ by the coordinate system contains the set of points with $\xi \leq \varepsilon+\delta$ and $\eta \leq \delta$. Let $\phi : (0,+\infty) \rightarrow (0,+\infty)$ be a function such that $\phi \in C^{\infty}(0,+\infty)$ and:
\begin{enumerate}
\item[(a)] $\phi(\xi) = 0$, if $\xi \leq \varepsilon$;
\item[(b)] $\phi(\xi) = \xi-\varepsilon$, if $\varepsilon+\delta \leq \xi$;
\item[(c)] $\phi(\xi)\leq \xi-\varepsilon$, for all $\xi \in (0,+\infty)$.
\end{enumerate}
With these choices we are able to redefine the projection for points classified as case $2$ on page $19$ of \cite{milnor}, i.e., $\varepsilon \leq \xi \leq \eta + \varepsilon$. Consider
\begin{equation*}
(t, u_1,\ldots,u_n) \mapsto (u_1,\ldots,u_{\lambda}, s_t u_{\lambda + 1},\ldots, s_t u_n),
\end{equation*}
where the number $s_t \in [0,1]$ is defined by
\begin{equation*}
s_t = t + (1-t) \sqrt{\frac{\phi(\xi)}{\eta}}.
\end{equation*}
This map is smooth for points in case $2$ because now the set $\phi(\xi)=\eta$ meets the boundary of $e^{\lambda}$ smoothly. By (c), we see that the image of each point at time $t=0$ is inside $M^{c-\varepsilon}\cup e^{\lambda}$. Using this new projection, the statement follows via the same program as in the proof of Theorem $3.2$ in \cite{milnor}.
\end{proof}

\bibliographystyle{amsbook}

\end{document}